\documentclass[12pt]{article}

\usepackage{amsmath}
\usepackage{amsthm,amscd,amsfonts}
\usepackage{amssymb, upref, color}
\usepackage{amsmath,amssymb,amsthm}
\usepackage{color}
\usepackage[colorlinks]{hyperref}
\usepackage[dvips]{graphicx}
\usepackage{epsf,epsfig,subfigure, verbatim}
\usepackage{latexsym,bm}

\numberwithin{equation}{section}

\theoremstyle{plain}
\newtheorem{exam}{Example}[section]
\newtheorem{theorem}[exam]{Theorem}
\newtheorem{lemma}[exam]{Lemma}
\newtheorem{remark}[exam]{Remark}

\newtheorem{proposition}[exam]{Proposition}
\newtheorem{definition}[exam]{Definition}
\newtheorem{corollary}[exam]{Corollary}

\begin{document}
\date{}

%\newpage

\title{ Solve the linear quaternion-valued differential equations having multiple eigenvalues
}
\author{ Kit Ian Kou$^{1}$\footnote{ Kit Ian Kou acknowledges financial support from the National Natural Science Foundation of China under Grant (No. 11401606), University of Macau (No. MYRG2015-00058-FST and No. MYRG099(Y1-L2)-FST13-KKI) and the Macao Science and Technology
Development Fund (No. FDCT/094/2011/A and No. FDCT/099/2012/A3). }
 \,\,\,\,\,  Wan-Kai Liu$^{1}$ \footnote{ Wankai Liu was supported by the National Natural
Science Foundation of China under Grant (No.11572288)}
\,\,\,\,\,\, Yong-Hui Xia$^{2,3}$
 \footnote{ Corresponding author. Email: xiadoc@163.com.  Yonghui Xia was supported by the National Natural
Science Foundation of China under Grant (No. 11671176 and No. 11271333), Natural
Science Foundation of Zhejiang Province under Grant (No. LY15A010007), Marie Curie Individual Fellowship within the European Community Framework Programme(MSCA-IF-2014-EF, ILDS - DLV-655209), the Scientific Research Funds of Huaqiao University and China Postdoctoral Science Foundation (No. 2014M562320). }
%\,\,\,\,\, Weinian Zhang$^{3}$\footnote{Weinian Zhang was supported by....}
\\
 {\small 1.Department of Mathematics, Faculty of Science and Technology, University of Macau, Macau}\\
{\small \em  kikou@umac.mo (K. I. Kou)  }
 \\
 {\small 2. School of Mathematical Sciences, Huaqiao University, 362021, Quanzhou, Fujian, China.}\\
{\small\em xiadoc@163.com  (Y-H.Xia)}\\
 {\small 3.Department of Mathematics, Zhejiang Normal University, Jinhua, 321004, China}
%{\small \em 3. School of Mathematics, Sicuan University, Chengdu, 610000, China}\\
%{\small E-mail: \,\,\,   matzwn@126.com (W.N. Zhang)    }
}

 \maketitle

\begin{center}
\begin{minipage}{140mm}
\begin{abstract}
The theory of two-dimensional linear quaternion-valued differential equations (QDEs) was recently established (see \cite{QDE1}). Some profound differences between QDEs and ODEs were observed.
 Also, an algorithm to evaluate the fundamental matrix by employing the eigenvalues and eigenvectors was presented in \cite{QDE1}. However, the fundamental matrix can be constructed providing that the eigenvalues are simple. If the linear system has multiple eigenvalues, how to construct the fundamental matrix? In particular, if the number of independent eigenvectors might be less than the dimension of the system. That is, the numbers of the eigenvectors are not enough to construct a fundamental matrix. How to find the ``missing solutions"?  The main purpose of this paper is to answer this question.

Furthermore, Caley determinant for Quaternion-valued matrix was adopted to proceed the theory of QDEs in \cite{QDE1}. One big disadvantage of Caley determinant is that it can be expanded along the different rows and columns. This may lead to different results due to non-commutativity of the quaternions. This approach is not convenient to be used. The novel definition of determinant for Quaternion-valued matrix based on permutation is introduced to analyze the theory.
This newly definition of determinant has great advantages compare to Calay determinant.

\end{abstract}



\end{minipage}
\end{center}

\section{\bf Introduction and Motivation}

%ͼƬǶÈë
%\begin{figure}
 % \centering
%  \includegraphics[scale=0.8]{a1.pdf}
%  \caption{Fig.4 An example of two limit cycle ([11]). }
%\end{figure}

Recently, Kou and Xia \cite{QDE1} established a framework for the basic theory of quaternion-valued differential equations (QDEs). QDEs have many applications in quantum mechanics (see e.g. \cite{Alder1,Alder2,Leo1,Leo2}), fluid mechanics (e.g. \cite{Gibbon1,Gibbon2,Gibbon3,RR1,RR2}), %Frenet frame in differential geometry (see e.g.\cite{Handson}), %kinematic modelling,
%attitude dynamics (see e.g. \cite{Wertz}), %Kalman filter design and spatial rigid body dynamics (\cite{Udwadia}),
etc.
However, there are few papers pursuing mathematical analysis for the QDEs. For examples,
Leo and Ducati \cite{Leo3} attempted to solve some special QDEs. Some results on the existence of periodic solutions were presented in \cite{Mawhin,W,W2}. Global analysis of a homogeneous QDE was given in Gasull et al. \cite{Zhang1} for $n=2,3$.
 A Bernoulli-type QDE was investigated by Zhang \cite{Zhang2}.
However, there is no systematic theory for QDEs.
Recently, Kou and Xia \cite{QDE1} presented a basic and theoretic framework for the two-dimensional linear homogeneous QDEs. They found four large differences between QDEs and ODEs due to the non-commutativity of the quaternion algebra. %The biggest difference is that the set of all the solutions to linear homogeneous QDEs is a right- free module.
 In \cite{QDE1}, the authors presented an algorithm to evaluate the fundamental matrix by employing the eigenvalue and eigenvectors. They provided a method and some examples to show how to construct the fundamental matrix when the eigenvalues are simple. However, it is possible to have multiple eigenvalues. How can we construct the fundamental matrix when the multiplicity of the eigenvalues is larger than one?  In particular, if the number of independent eigenvectors might be smaller than the dimensionality of the system. That is, the numbers of the eigenvectors is not enough to construct a fundamental matrix. In this paper, one of the main tasks is to find  the ``missing solutions". We will devote ourselves to answer this question.

On the other hand, Cayley determinant for quaternion valued matrix \cite{Cayley} was adopted in \cite{QDE1} which depends on the expansion of $i$-th row and $j$-th column of the quaternion-valued matrix.  Different expansions of the quaternionic determinant can lead to different results.
For example, the results of expanding along the column and expanding along the row are different.
  Owing to the non-commutativity of the quaternion algebra, the results are different due to different expansions. The determinant of $n$-order matrix is more complicate. Therefor Cayley determinant is not convenient for the quaternion valued matrix. We will adopt another definition (see eg. Chen \cite{Chen1}) to analyze our results in this paper.
This definition of determinant which is based on permutation has great advantage compared to Calay determinant (see next section in detail). Due to the newly definition of determinant, the computation of the determinant is different. In particular, the proof of Liouville formula is more complicated.

\section{\bf Quaternion algebra}
For the quaternion algebra, we adopt the operators and notations in \cite{QDE1}. To avoid repeating here, we omit some standard definitions (e.g. conjugate, norm) which has been introduced in \cite{QDE1}.

As pointed out in the introduction, Caley determinant can be expanded along the $j$-th column or the $i$-th row.
 Owning to the non-commutativity of the quaternion algebra, the results of the determinant are not same due to different expansions. Thus, it is not convenient to apply this definition to the quaternion valued matrix.
 In 1991, Chen \cite{Chen1} gave us a ``direct" definition by specifying a certain ordering of the factors in the $n!$ terms in the sum.
 %And the definition  introduced below, of a double determinant, is the one that gives us many properties resembling ordinary ones \cite{Chen1,Chen2,Chen3}.
 In this sense, the determinant has a unique result. So, in this article, we will study the $n$ dimensional linear quaternionic-valued ordinary differential equations based on this definition. Now we are in a position to introduce this definition of determinant.

Let $\mathbb{H}^{n\times m}$ denote the set of all matrices $A=(a_{i,j})_{n\times m}$, where $a_{i,j}$ are quaternions. For any $A\in \mathbb{H}^{n\times n}$, the determinant based on permutation is defined as follows (see e.g. \cite{Chen1}).
\[
\det_{P} A=|A|_{P}= \det_{P}\left ( \begin{array}{ll}   {a} _{11}  & {a} _{12}  \cdots  {a} _{1n}
 \\
 {a} _{21}  & {a} _{22}  \cdots  {a} _{2n}
 \\
 \,\,\,\vdots &\,\,\,\vdots\,\,\,\ddots \,\,\,\vdots
 \\
{a} _{n1}  & {a} _{n2}  \cdots  {a} _{nn}
\end{array}
\right )
\]
\begin{equation}
\equiv\sum_{\sigma\in S_{n}}\varepsilon(\sigma)a_{n_{1}i_{2}}a_{i_{2}i_{3}}\cdots a_{i_{s-1}i_{s}}a_{i_{s}n_{1}}a_{n_{2}j_{2}}\cdots a_{j_{t}n_{2}}\cdots a_{n_{r}k_{2}}\cdots a_{k_{l}n_{r}},
\label{jzdy}
\end{equation}
where $S_{n}$ is the symmetric group on $n$ letters, and the disjoint cycle decomposition of $\sigma\in S_{n}$ is written in the normal form:
\[
\sigma=(n_{1}i_{2}i_{3}\cdots i_{s})(n_{2}j_{2}j_{3}\cdots j_{t})\cdots (n_{r}k_{2}k_{3}\cdots k_{l}),
\]
\[
n_{1}>i_{2},i_{3},\cdots, i_{s},n_{2}>j_{2},j_{3},\cdots, j_{t},\cdots , n_{r}>k_{2},k_{3},\cdots, k_{l},
\]
\[
n=n_{1}>n_{2}>\cdots >n_{r}\geq 1,
\]
and
\[
\varepsilon(\sigma)=(-1)^{(s-1)+(t-1)+\cdots+(l-1)}=(-1)^{n-r}.
\]
Different from the Caley determinant $\det$, we denote this kind of determinant based on permutation by $\det\limits_{P}$.
Notice that if all $a_{ij}$ commute with each other, the definition of $\det A$ is the same as that an ordinary determinant (Caley  determinant).

For convenience and explicitness, %if cycle factor $\sigma_{0}=(i_{1}i_{2}\cdots i_{s})$, and
we denote
\[
a_{i_{1}i_{2}}a_{i_{2}i_{3}}\cdots a_{i_{s-1}i_{s}}=<i_{1}i_{2}\cdots i_{s}>,
\]
and
\[
a_{i_{1}i_{2}}a_{i_{2}i_{3}}\cdots a_{i_{s-1}i_{s}}a_{i_{s}i_{1}}=<i_{1}i_{2}\cdots i_{s}i_{1}>.
\]
Then expression (\ref {jzdy}) is simplified into
\[
\det\limits_{P} A=\sum_{\sigma\in S_{n},\sigma=\sigma_{1}\sigma_{2}\cdots\sigma_{r}}\varepsilon(\sigma)<\sigma_{1}><\sigma_{2}>\cdots<\sigma_{r}>=\sum_{\sigma\in S_{n}}\varepsilon(\sigma)<\sigma>.
\]
In particular, for $n=2$,
\[
 \sigma_{1}=(2)(1),\sigma_{2}=(21)\in S_{2},
\]
and
\[
 \varepsilon(\sigma_{1})=(-1)^{(1-1)+(1-1)}=1,\varepsilon(\sigma_{2})=(-1)^{(2-1)}=-1.
\]
Consequently, we have
\[
\det\limits_{P} \left ( \begin{array}{cc}
a_{11} & a_{12}
\\
a_{21} & a_{22}
\end{array}
\right )=\varepsilon(\sigma_{1})a_{22}a_{11}+\varepsilon(\sigma_{2})a_{21}a_{12}=a_{22}a_{11}-a_{21}a_{12}.
\]
For $n=3$, we have
\[
 \sigma_{1}=(3)(2)(1),\sigma_{2}=(3)(21),\sigma_{3}=(312),
\]
\[
 \sigma_{4}=(31)(2),\sigma_{5}=(321),\sigma_{6}=(32)(1)\in S_{3},
\]
and
\[
 \varepsilon(\sigma_{1})=(-1)^{(3-1)}=1,\varepsilon(\sigma_{2})=(-1)^{(2-1)}=-1,
 \varepsilon(\sigma_{3})=(-1)^{(1-1)}=1
\]
\[
 \varepsilon(\sigma_{4})=(-1)^{(2-1)}=-1,\varepsilon(\sigma_{5})=(-1)^{(1-1)}=-1,
 \varepsilon(\sigma_{6})=(-1)^{(2-1)}=1.
\]
So we have
\[
\det_P\left (\begin{array}{ccc}
a_{11} & a_{12} & a_{13}
\\
a_{21} & a_{22} & a_{23}
\\
a_{31}& a_{32} & a_{33}
\end{array}
\right )
=
%\begin{array{lll}
%&=&
%\varepsilon(\sigma_{1})a_{33}a_{22}a_{11}+\varepsilon(\sigma_{2})a_{33}a_{21}a_{12}+\varepsilon(\sigma_{3})
%a_{31}a_{12}a_{23}+\varepsilon(\sigma_{4})a_{31}a_{13}a_{22}+\varepsilon(\sigma_{5})a_{32}a_{21}a_{13}+
%\varepsilon(\sigma_{6})a_{32}a_{23}a_{11}
%\\
%&=&
 a_{33}a_{22}a_{11}-a_{33}a_{21}a_{12}+a_{31}a_{12}a_{23}
-a_{31}a_{13}a_{22}+a_{32}a_{21}a_{13}-a_{32}a_{23}a_{11}.
%\end{array}
\]
%\[
%\left | \begin{array}{ccc} 1 & \boldsymbol{i} & \boldsymbol{j}
%\\
%-\boldsymbol{i} & 2 & \boldsymbol{j}-\boldsymbol{k}
%\\
%-\boldsymbol{j} & -\boldsymbol{j}+\boldsymbol{k} & 3
%\end{array}
%\right %|=6-3(-\boldsymbol{i})\boldsymbol{i}-(-\boldsymbol{j}+\boldsymbol{k})(\boldsymbol{j}-\boldsymbol{k})-2(\boldsymbol{j})\boldsymbol{j}+(-\boldsymbol{j}+\boldsymbol{k})(-\boldsymbol{i})\boldsymbol{j}+(-\boldsymbol{j})\boldsymbol{i}(\boldsymbol{j}-\boldsymbol{k})=1
%\]

For any $A\in \mathbb{H}^{n\times m}$, $\rm{ddet}_P A \equiv\det\limits_{P}(A^{+}A)$ is called double determinant of $A$. Since $A^{+}A$ is a Hermitian matrix, $\rm{ddet}_PA$ is always a real number (it can be proved that $\rm{ddet}_PA \geq 0$).
%For above $3\times3$ determinant, we can get the value of determinant is $2\boldsymbol{i}-1$ (expanding along the first row) and $-2\boldsymbol{i}-1$ (expanding along the first column), which will cause us a lot of trouble in the course of research next. this is an important reason for we replacing the Cayley determinant by above definition.

Let $\psi : \mathbb{R} \rightarrow \mathbb{H}$ be a quaternion-valued function defined on $\mathbb{R}$. We
denote the set of such quaternion-valued functions by $\mathbb{H} \otimes \mathbb{R}$. Then $n$ dimensional quaternionic functions of real variable, $\mathbb{H}^{n}\otimes\mathbb{R}= \{\Psi(t) \big| \Psi(t) =(\psi_{1}(t),\psi_{2}(t),\cdots,\psi_{n}(t))^{T} \}$. And the derivative, integral and norms of $n$ dimensional quaternionic functions with respect to the real variable $t$ are well defined, one can refer to \cite{QDE1}. Moverover, we adopt these notations from \cite{QDE1}.

\section{Wronskian and Structure of General Solution to QDEs}

Consider the $n$-dimensional linear QDEs as follows.
\begin{equation}
\dot \Psi (t) = A(t)\Psi(t)\,\,\,\, or \,\,\,\,\,\,
\left ( \begin{array}{l} \dot {\psi_{1}} (t)
\\
\dot {\psi_{2}} (t)
\\
\cdots
\\
\dot{\psi_{n}} (t)
\end{array}
\right )= \left ( \begin{array}{ll}   {a} _{11}(t)  & {a} _{12}(t)  \cdots  {a} _{1n}(t)
 \\
 {a} _{21}(t)  & {a} _{22}(t)  \cdots  {a} _{2n}(t)
 \\
               & \cdots\cdots
 \\
{a} _{n1}(t)  & {a} _{n2}(t)  \cdots  {a} _{nn}(t)
\end{array}
\right ) \left ( \begin{array}{l}   {\psi_{1}} (t)
\\
{\psi_{2}} (t)
\\
\cdots
\\
{\psi_{n}} (t)
\end{array}
\right ),
\label{2.1}
\end{equation}
where $\Psi (t) \in \mathbb{H}^n \otimes \mathbb{R}$, $A(t)  \in   \mathbb{H}^{n\times n} \otimes \mathbb{R}$ is continuous on the interval $[a,b]$.
 Similar discussion to Theorem 3.1 in \cite{QDE1}, we have

 \begin{theorem}\label{Th3.1}
System $(\ref{2.1})$  has exactly one solution satisfying the following initial value problem
\[
\Psi (t_0) =  \xi, \,\,\,\,\,\,  \xi  \in \mathbb{H}^n .
\label{2.2}
\]
 \end{theorem}

Now we introduce some definitions on abelian groups, rings, modules, submodules, direct sum in the abstract algebraic theory. % (e.g. \cite{abstract,abstract2}).
 Due to the length of this paper, we omit the definitions. For these detailed definitions, one can refer to \cite{QDE1}.
We also adopt the notations from \cite{QDE1}. %As defined in \cite{QDE1}, the elements $x_{1}, \cdots, x_{k}\in \mathbb{A}_{\cal{R}}^{R}$ are said to be right independent if
%\[
%x_1 r_1+ \cdots+x_k r_k=0,r_i\in\mathcal R~\text{implies that}~r_1=\cdots=r_k=0.
%\]
%that is, the only linear combination that vanishes is the trivial one with every coefficient zero.
%For the sake of convenience, we call it {\em right} independence.
%A subset $\{x_{1}, \cdots, x_{k}\}$ is called a basis of a right module $\mathbb{A}_{\cal{R}}^{R}$ if it is right independent and generates  $\mathbb{A}_{\cal{R}}^{R}$. A module that has a finite basis is called a {\bf free module}.

%Direct differentiation, one can verify that

%\begin{theorem}\label{Th4.1}
%$($superposition theorem$)$ If $u(t)$ and $v(t)$ are two solutions of Eq.$($\ref{2.1}$)$, then the linear combination $u(t)\alpha+v(t)\beta$ is also a solution of Eq.$($\ref{2.1}$)$, where $\alpha,\beta\in \mathbb{H}$.
%\end{theorem}

%By Theorem \ref{Th4.1}, it is possible that all solutions of Eq.(\ref{2.1}) can generate the right $\mathbb{H}$-module.
We claim that the set of all the solutions to Eq.(\ref{2.1}) is the right $\mathbb{H}$-module. To prove this, firstly we try to find a basis of this right $\mathbb{H}$-module. Now we should introduce the concept of independence and  dependence for the vector functions $x_1(t),x_2(t),\cdots,x_n(t)$.

\begin{definition}\label{De4.1}
For $n$ quaternion-valued vector functions $x_1(t),x_2(t),\cdots,x_n(t)$, each $x_i(t)\in\mathbb{H}^{n}\otimes \mathbb{R}$ defined on the real interval $I$, if
\[
x_1(t) r_1+ \cdots+x_n(t) r_n=0,r_i\in\mathbb{H}~\text{implies that}~r_1=\cdots=r_k=0, t\in I
\]
then $x_1(t),x_2(t),\cdots,x_n(t)$ is said to be independent. Otherwise, $x_1(t),x_2(t),\cdots,x_n(t)$ is said to be dependent.
\end{definition}

\begin{definition}\label{De4.2}
Let $x_1(t),x_2(t), \cdots , x_n(t)$ are n solutions of Eq.$($\ref{2.1}$)$. Let
\[
M(t)
=(x_1(t),x_2(t), \cdots , x_n(t)).
%\left ( \begin{array}{ll}
%{x} _{11}(t)  & {x} _{12}(t) \cdots {x} _{1n}(t)
%\\
%{x} _{21}(t)  & {x} _{22}(t) \cdots {x} _{2n}(t)
%\\
%              & \cdots\cdots
%\\
% {x} _{n1}(t)  & {x} _{n2}(t) \cdots {x} _{nn}(t)
%\end{array}
%\right ).
\]
The {\em Wronskian} of QDEs is defined by
\[
 W_{QDE}^P(t) = \rm{ddet}_P M(t): =  \det_P \big( M^+(t) M(t)\big),
\]
where $M^{+}$ is the conjugate transpose of $M(t)$, namely
\[
M^{+}(t)
=
\left ( \begin{array}{ll}
 \overline {x} _{11}(t)  &  \overline {x} _{21}(t) \cdots \overline {x} _{n1}(t)
\\
\overline {x} _{12}(t)  & \overline {x} _{22}(t) \cdots \overline {x} _{n2}(t)
 \\
 & \cdots\cdots
 \\
\overline {x} _{1n}(t)  & \overline {x} _{2n}(t) \cdots \overline {x} _{nn}(t)
\end{array}
\right ).
\]
\end{definition}

\begin{remark}\label{Re4.2}
As pointed out in \cite{QDE1}, the standard {\em Wronskian} of ODEs is not valid for QDEs. So we define {\em Wronskian} of QDEs by $\det\limits_P \big( M^+(t) M(t)\big)$.
For quaternion matrix, it should be noted that
\[
\det_P\big(A(t) B(t)\big)\neq \det_P A(t) \cdot \det_P   B(t)  .
\]
But, in \cite{Chen2} $($Theorem 5$)$, he proved that
\[
{\rm {ddet}_P} \big(A(t) B(t)\big)= {\rm {ddet}_P} A(t) \cdot {\rm {ddet}_P}   B(t)  .
\]
\end{remark}

\begin{lemma}\label{Th4.2}
If $x_1(t), x_2(t), \cdots , x_n(t) $ are right dependent on $I$, then $W_{QDE}^P(t)=0$.
\end{lemma}

The proof is standard (similar to ODEs). We omit the proof.

%{\bf Proof.} If $x_1(t), x_2(t), \cdots , x_n(t) $ are right dependent on $I$, then there exist $r_{1},r_{2}, \cdots , r_{n}\in\mathbb{H}$ (not all zero) such that
%\[
%x_1(t) r_1 + x_2(t) r_2 + \cdots + x_n(t) r_n =0,\,\,\,\,\ t\in I,
%\]
%or
%\begin{equation}
%\left( x_1(t), x_2(t), \cdots,x_n(t)\right)\left ( \begin{array}{ll}
% r_1
%\\
%r_2
% \\
% \cdots
% \\
%r_n
%\end{array}
%\right )  =0,\,\,\,\,\ t\in I.
%\label{2.4}
%\end{equation}
%In Eq.(\ref{2.4}), $\left( x_1(t), x_2(t), \cdots,x_n(t)\right)$ can be seen as the coefficient matrix, $r_{1},r_{2},\cdots,r_{n}$ can be seen as unknown numbers. Thus, Eq.(\ref{2.4}) can be seen as homogeneous system of right linear equations $M(t)r=0$. By way of contradiction, if there exists $t_{0}\in I $ such that $W_{QDE}^P(t_0)\neq0$, by Theorem 2.3 \cite{Chen1}, then exists a unique solution $r_{1}=r_{2}=\cdots=r_{n}=0$. This implies that $x_1(t), x_2(t), \cdots , x_n(t) $ are right dependent on $I$, which is a contradiction. Therefore, $W_{QDE}^P(t)=0$, $t\in I$.

We need a lemma from Theorem 8 in \cite{Chen2}.

\begin{lemma}\label{Lemma8}
Let $A_{n\times m}=(\alpha_1, \alpha_2,\cdots, \alpha_m)^T\in \mathbb{H}^{n\times m}$, where $\alpha_1, \alpha_2,\cdots, \alpha_m$ are $m$ column quaternionic vectors. Then $\alpha_1, \alpha_2,\cdots, \alpha_m$ are right independent if and only if $ddet A_{n\times m}\neq0$.
\end{lemma}

In particular, for $n=m$, then $\alpha_1, \alpha_2,\cdots, \alpha_n$ are right independent if and only if $ddet A_{n\times n}\neq0$.

\begin{lemma}\label{Th4.3}
 If $x_1(t), x_2(t), \cdots , x_n(t) $ are $n$ right independent solutions of Eq.$($\ref{2.1}$)$, then  $W_{QDE}^P(t)\neq0$ on $I$.
 \end{lemma}

 The proof is standard (similar to ODEs). We omit the proof. Now, we present Liouville formula as follows.

\begin{theorem}\label{Th4.4}
%$($Liouville formula$)$ The {\em Wronskian}  $W_{QDE}^P(t)$ of Eq. $($\ref{2.1}$)$ satisfies the following quaternionic Liouville formula.
\[
W_{QDE}^P(t)=exp\left(\int_{t_{0}}^{t}[trA(s)+trA^{+}(s)]ds\right)W_{QDE}^P(t_{0}),
\]
or
\[
W_{QDE}^P(t)=exp\left(\int_{t_{0}}^{t}2\Re(trA(t))ds\right)W_{QDE}^P(t_{0}),
\]
where $trA(t)$ is the trace of the coefficient matrix $A(t)$, i.e. $trA(t)= \sum\limits^{n}_{i=1} a_{ii}(t)$.% Moreover, if $W_{QDE}^P(t)= 0$ at some $t_{0}$ in I then $W_{QDE}^P(t)= 0$ on I.
\end{theorem}

{\bf Proof.} By definition of determinant
\[
\frac{d}{dt}W_{QDE}^P(t)
=\frac{d}{dt}\left| \left ( \begin{array}{ll}
\overline{x} _{11}(t) & \overline{x} _{21}(t) \,\, \cdots \,\, \overline{x} _{n1}(t)
\\
\overline{x} _{12}(t) & \overline{x} _{22}(t) \,\, \cdots \,\, \overline{x} _{n2}(t)
\\
\,\,\,\,\,\,\vdots \,\,\, & \,\,\,\,\,\, \vdots  \,\,\,\,\,\,\,\,\,\ddots \,\,\,\,\,\,\,\,\, \vdots
\\
\overline{x} _{1n}(t) & \overline{x} _{2n}(t)\, \, \cdots \,\, \overline{x} _{nn}(t)
\end{array}
\right)
\left ( \begin{array}{ll}
{x} _{11}(t) & {x} _{12}(t) \,\, \cdots \,\, {x} _{1n}(t)
\\
{x} _{21}(t) & {x} _{22}(t) \,\, \cdots \,\, {x} _{2n}(t)
 \\
 \,\,\,\,\,\,\vdots \,\,\, & \,\,\,\,\,\, \vdots  \,\,\,\,\,\,\,\,\,\ddots \,\,\,\,\,\,\,\,\, \vdots
 \\
{x} _{n1}(t) & {x} _{n2}(t) \,\, \cdots \,\, {x} _{nn}(t)
\end{array}
\right)
\right|_{P}
\]
\[
=\frac{d}{dt}\left | \begin{array}{ll}
\sum\limits_{i=1}^{n} \overline {x} _{i1}(t){x} _{i1}(t) & \sum\limits_{i=1}^{n} \overline {x} _{i1}(t){x} _{i2}(t) \,\, \cdots \,\, \sum\limits_{i=1}^{n} \overline {x} _{i1}(t){x} _{in}(t)
\\
\sum\limits_{i=1}^{n} \overline {x} _{i2}(t){x} _{i1}(t) & \sum\limits_{i=1}^{n} \overline {x} _{i2}(t){x} _{i2}(t) \,\, \cdots \,\, \sum\limits_{i=1}^{n} \overline {x} _{i2}(t){x} _{in}(t)
\\
\,\,\,\,\,\,\,\,\,\,\,\,\,\,\,\,\, \vdots \,\,\,\,\,\,\,\,\,\,\,\,\,\,\,\, & \,\,\,\,\,\,\,\,\,\,\,\,\,\,\,\,\,\, \vdots
\,\,\,\,\,\,\, \,\,\,\,\,\,\,\,\,\,\, \ddots \,\,\,\,\,\,\, \,\,\,\,\,\,\,\,\,\,\,\,\,\, \vdots
\\
\sum\limits_{i=1}^{n} \overline {x} _{in}(t){x} _{i1}(t) & \sum\limits_{i=1}^{n} \overline {x} _{in}(t){x} _{i2}(t) \,\, \cdots \,\, \sum\limits_{i=1}^{n} \overline {x} _{in}(t){x} _{in}(t)
\end{array}
\right|_{P}
\]
\[
=\left | \begin{array}{ll}
\frac{d}{dt}\sum\limits_{i=1}^{n} \overline {x} _{i1}(t){x} _{i1}(t)  \,\, \cdots \,\, \frac{d}{dt}\sum\limits_{i=1}^{n} \overline {x} _{i1}(t){x} _{in}(t)
\\
\,\,\,\,\,\,\sum\limits_{i=1}^{n} \overline {x} _{i2}(t){x} _{i1}(t)  \,\, \cdots \,\,\,\,\,\,\,\, \sum\limits_{i=1}^{n} \overline {x} _{i2}(t){x} _{in}(t)
\\
\,\,\,\,\,\,\,\,\,\,\,\,\,\,\,\,\,\,\,\,\,\,\, \vdots \,\,\,\,\,\,\,\, \,\,\,\,\,\,\,\,\,\,\, \ddots \,\,\,\,\,\,\,\,\,\,\, \,\,\,\,\,\,\,\,\,\,\,\,\,\, \vdots
\\
\,\,\,\,\,\,\sum\limits_{i=1}^{n} \overline {x} _{in}(t){x} _{i1}(t)  \,\, \cdots \,\,\,\,\,\,\,\, \sum\limits_{i=1}^{n} \overline {x} _{in}(t){x} _{in}(t)
\end{array}
\right|_{P}
+\left | \begin{array}{ll}
\,\,\,\,\,\,\sum\limits_{i=1}^{n} \overline {x} _{i1}(t){x} _{i1}(t)  \,\, \cdots \,\,\,\,\,\,\,\, \sum\limits_{i=1}^{n} \overline {x} _{i1}(t){x} _{in}(t)
\\
\frac{d}{dt}\sum\limits_{i=1}^{n} \overline {x} _{i2}(t){x} _{i1}(t)  \,\, \cdots \,\, \frac{d}{dt}\sum\limits_{i=1}^{n} \overline {x} _{i2}(t){x} _{in}(t)
\\
\,\,\,\,\,\,\,\,\,\,\,\,\,\,\,\,\,\,\,\,\,\,\, \vdots \,\,\,\,\,\,\,\, \,\,\,\,\,\,\,\,\,\,\, \ddots \,\,\,\,\,\,\,\,\,\,\, \,\,\,\,\,\,\,\,\,\,\,\,\,\, \vdots
\\
\,\,\,\,\,\,\sum\limits_{i=1}^{n} \overline {x} _{in}(t){x} _{i1}(t)  \,\, \cdots \,\,\,\,\,\,\,\, \sum\limits_{i=1}^{n} \overline {x} _{in}(t){x} _{in}(t)
\end{array}
\right|_{P}
\]
\begin{equation}
+\cdots+ \left | \begin{array}{ll}
\,\,\,\,\,\,\sum\limits_{i=1}^{n} \overline {x} _{i1}(t){x} _{i1}(t)  \,\, \cdots \,\,\,\,\,\,\,\, \sum\limits_{i=1}^{n} \overline {x} _{i1}(t){x} _{in}(t)
\\
\,\,\,\,\,\,\sum\limits_{i=1}^{n} \overline {x} _{i2}(t){x} _{i1}(t)  \,\, \cdots \,\,\,\,\,\,\,\, \sum\limits_{i=1}^{n} \overline {x} _{i2}(t){x} _{in}(t)
\\
\,\,\,\,\,\,\,\,\,\,\,\,\,\,\,\,\,\,\,\,\,\,\, \vdots \,\,\,\,\,\,\,\, \,\,\,\,\,\,\,\,\,\,\, \ddots \,\,\,\,\,\,\,\,\,\,\, \,\,\,\,\,\,\,\,\,\,\,\,\,\, \vdots
\\
\frac{d}{dt}\sum\limits_{i=1}^{n} \overline {x} _{in}(t){x} _{i1}(t)  \,\, \cdots \,\, \frac{d}{dt}\sum\limits_{i=1}^{n} \overline {x} _{in}(t){x} _{in}(t)
\end{array}
\right|_{P}
\label{0.1}
\end{equation}
Since $x_{1}(t),x_{2}(t),\cdots,x_{k}(t)$, $k=1,2 \cdots n$ are $n$ solutions of the Eq.(\ref{2.1}), we have
\[
\left | \begin{array}{ll}
\frac{d}{dt}\sum\limits_{i=1}^{n} \overline {x} _{i1}(t){x} _{i1}(t)  \,\, \cdots \,\, \frac{d}{dt}\sum\limits_{i=1}^{n} \overline {x} _{i1}(t){x} _{in}(t)
\\
\,\,\,\,\,\,\sum\limits_{i=1}^{n} \overline {x} _{i2}(t){x} _{i1}(t)  \,\, \cdots \,\,\,\,\,\,\,\, \sum\limits_{i=1}^{n} \overline {x} _{i2}(t){x} _{in}(t)
\\
\,\,\,\,\,\,\,\,\,\,\,\,\,\,\,\,\,\,\,\,\,\,\, \vdots \,\,\,\,\,\,\,\, \,\,\,\,\,\,\,\,\,\,\, \ddots \,\,\,\,\,\,\,\,\,\,\, \,\,\,\,\,\,\,\,\,\,\,\,\,\, \vdots
\\
\,\,\,\,\,\,\sum\limits_{i=1}^{n} \overline {x} _{in}(t){x} _{i1}(t)  \,\, \cdots \,\,\,\,\,\,\,\, \sum\limits_{i=1}^{n} \overline {x} _{in}(t){x} _{in}(t)
\end{array}
\right|_{P}=
\]
\[
\left | \begin{array}{ll}
\sum\limits_{i=1}^{n}\sum\limits_{j=1}^{n}(\overline {x} _{j1}(t)\overline {a} _{ij}(t){x} _{i1}(t)+\overline{x} _{i1}(t) {a} _{ij}(t){x} _{j1}(t))  \cdots  \sum\limits_{i=1}^{n}\sum\limits_{j=1}^{n}(\overline {x} _{j1}(t)\overline {a} _{ij}(t){x} _{in}(t)+\overline{x} _{i1}(t) {a} _{ij}(t){x} _{jn}(t))
\\
\,\,\,\,\,\,\,\,\,\,\,\,\,\,\,\,\,\,\,\,\,\,\,\,\,\,\,\,\,\,\,\,\,\,\,\,\,\,\,\,\,\,\,\sum\limits_{i=1}^{n} \overline {x} _{i2}(t){x} _{i1}(t)  \,\,\,\,\,\,\,\,\,\,\,\,\,\,\,\,\,\,\,\,\,\,\,\,\,\,\,\,\,\,\,\, \cdots \,\,\,\,\,\,\,\,\,\,\,\,\,\,\,\,\,\,\,\,\,
\,\,\,\,\,\,\,\,\,\,\,\,\,\,\,\,\,\,\,\,\,\,\,\sum\limits_{i=1}^{n} \overline {x} _{i2}(t){x} _{in}(t)
\\
\,\,\,\,\,\,\,\,\,\,\,\,\,\,\,\,\,\,\,\,\,\,\,\,\,\,\,\,\,\,\,\,\,\,\,\,\,\,\,\,\,\,\,\,\,\,\,\,\,\,\,\,\,\,\,\,\,\,\,\, \vdots \,\,\,\,\,\,\,\,\,\,\,\,\,\,\,\,\,\,\,\,\,\,\,\,\,\,\,\,\,\,\,\,\,\, \,\,\,\,\,\,\,\,\,\,\,\,\,\,\, \ddots \,\,\,\,\,\,\,\,\,\,\,\,\,\,\,\,\,\,\,\,\,\,\,\,\,\,\,\,\,\,\,\,\,\,\,\,\,\,\,\,\,\,\,\,\,\,\,\, \,\,\,\,\,\,\,\,\,\,\,\,\,\, \vdots
\\
\,\,\,\,\,\,\,\,\,\,\,\,\,\,\,\,\,\,\,\,\,\,\,\,\,\,\,\,\,\,\,\,\,\,\,\,\,\,\,\,\,\,\sum\limits_{i=1}^{n} \overline {x} _{in}(t){x} _{i1}(t)  \,\, \,\,\,\,\,\,\,\,\,\,\,\,\,\,\,\,\,\,\,\,\,\,\,\,\,\,\,\,\,\, \cdots \,\,\,\,\,\,\,\,\,\,\,\,\,\,\,\,\,\,\,\,\,
\,\,\,\,\,\,\,\,\,\,\,\,\,\,\,\,\,\,\,\,\,\,\, \sum\limits_{i=1}^{n} \overline {x} _{in}(t){x} _{in}(t)
\end{array}
\right|_{P}
\]
\[
=\left | \begin{array}{ll}
2\Re a_{ii} \sum\limits_{i=1}^{n} \overline {x} _{i1}(t){x} _{i1}(t) & 2\Re a_{ii}\sum\limits_{i=1}^{n} \overline {x} _{i1}(t){x} _{i2}(t) \,\, \cdots \,\, 2\Re a_{ii}\sum\limits_{i=1}^{n} \overline {x} _{i1}(t){x} _{in}(t)
\\
\,\,\,\,\,\sum\limits_{i=1}^{n} \overline {x} _{i2}(t){x} _{i1}(t) & \,\,\,\,\,\sum\limits_{i=1}^{n} \overline {x} _{i2}(t){x} _{i2}(t) \,\, \cdots \,\, \,\,\,\,\,\sum\limits_{i=1}^{n} \overline {x} _{i2}(t){x} _{in}(t)
\\
\,\,\,\,\,\,\,\,\,\,\,\,\,\,\,\,\,\,\,\,\,\, \vdots \,\,\,\,\,\,\,\,\,\,\,\,\,\,\,\, & \,\,\,\,\,\,\,\,\,\,\,\,\,\,\,\,\,\,\,\,\,\,\, \vdots
\,\,\,\,\,\,\, \,\,\,\,\,\,\,\,\,\,\, \ddots \,\,\,\,\,\,\,\,\, \,\,\,\,\,\,\,\,\,\,\,\,\,\, \vdots
\\
\,\,\,\,\,\sum\limits_{i=1}^{n} \overline {x} _{in}(t){x} _{i1}(t) & \,\,\,\,\,\sum\limits_{i=1}^{n} \overline {x} _{in}(t){x} _{i2}(t) \,\, \cdots \,\, \,\,\,\,\,\sum\limits_{i=1}^{n} \overline {x} _{in}(t){x} _{in}(t)
\end{array}
\right|_{P}
\]
\[
+\sum\limits_{k=1}^{n}\sum\limits_{\substack{j=1\\ j \neq k}}^{n}
\left | \begin{array}{ll}
\overline {x} _{j1}(t)\overline {a} _{kj}(t){x} _{k1}(t)+\overline{x} _{k1}(t) {a} _{kj}(t){x} _{j1}(t)  \cdots  \overline {x} _{j1}(t)\overline {a} _{kj}(t){x} _{in}(t)+\overline{x} _{k1}(t) {a} _{kj}(t){x} _{jn}(t)
\\
\,\,\,\,\,\,\,\,\,\,\,\,\,\,\,\,\,\,\,\,\,\,\,\,\,\,\,\,\,\sum\limits_{i=1}^{n} \overline {x} _{i2}(t){x} _{i1}(t)  \,\,\,\,\,\,\,\,\,\,\,\,\,\,\,\,\,\,\,\,\,\,\,\,\,\,\,\,\, \cdots \,\,\,\,\,\,\,
\,\,\,\,\,\,\,\,\,\,\,\,\,\,\,\,\,\,\,\,\,\sum\limits_{i=1}^{n} \overline {x} _{i2}(t){x} _{in}(t)
\\
\,\,\,\,\,\,\,\,\,\,\,\,\,\,\,\,\,\,\,\,\,\,\,\,\,\,\,\,\,\,\,\,\,\,\,\,\,\,\,\,\,\,\,\,\,\, \vdots \,\,\,\,\,\,\,\,\,\,\,\,\,\,\,\,\,\,\,\,\,\,\,\,\,\,\,\,\,\,\,\,\,\,\,\,\,\,\,\,\,\,\,\,\,  \ddots \,\,\,\,\,\,\,\,\,\,\,\,\,\,\,\,\,\,\,\,\,\,\,\,\,\,\,\,\,\,\,\,\,\,\,\,\,\,\,\,\,\,\,\,\,\, \, \vdots
\\
\,\,\,\,\,\,\,\,\,\,\,\,\,\,\,\,\,\,\,\,\,\,\,\,\,\,\,\,\sum\limits_{i=1}^{n} \overline {x} _{in}(t){x} _{i1}(t)  \,\, \,\,\,\,\,\,\,\,\,\,\,\,\,\,\,\,\,\,\,\,\,\,\,\,\,\,\,\cdots \,\,\,\,\,\,\,\,
\,\,\,\,\,\,\,\,\,\,\,\,\,\,\,\,\,\,\,\, \sum\limits_{i=1}^{n} \overline {x} _{in}(t){x} _{in}(t)
\end{array}
\right|_{P}
\]
$:=A_{1}(t)+\sum\limits_{k=1}^{n}\sum\limits_{\substack{j=1\\ j \neq k}}^{n}B_{kj}^{1}(t)$.\\
%where $b_{i}=a_{ii}+\overline a_{ii}$.
 Similarly, the expansion of the $m$ determinant of the Eq.(\ref{0.1}) is $A_{m}(t)+\sum\limits_{k=1}^{n}\sum\limits_{\substack{j=1\\ j \neq k}}^{n}B_{kj}^{m}(t)$, ($m=1,2\cdots n$). So we get
\begin{equation}
\frac{d}{dt}W_{QDE}^P(t)=\sum\limits_{m=1}^{n}A_{m}(t)+\sum\limits_{m=1}^{n}\sum\limits_{k=1}^{n}\sum\limits_{\substack{j=1\\ j \neq k}}^{n}B_{kj}^{m}(t)=\sum\limits_{m=1}^{n}A_{m}(t)+
\sum\limits_{k=1}^{n}\sum\limits_{\substack{j=1\\j \neq k}}^{n}\sum\limits_{m=1}^{n}B_{kj}^{m}(t).
\label{0.2}
\end{equation}

(i) Considering $W_{QDE}^P(t)$. For every permutation $\sigma\in S_{n}$,
\[
\sigma=(n_{1}i_{2}\cdots i_{s})\cdots(n_{p} j_{2} \cdots j_{q})\cdots(n_{r} k_{2} \cdots k_{t})=\sigma_{1}\cdots\sigma_{p}\cdots\sigma_{r},
\]
we have
\begin{equation}
\begin{array}{lll}
<\sigma>
&=&
(\sum\limits_{i=1}^{n}\overline {x}_{in_{1}}x_{ii_{2}})\cdots(\sum\limits_{i=1}^{n}\overline {x}_{ii_{s}}x_{in_{1}})\cdots
(\sum\limits_{i=1}^{n}\overline {x}_{in_{p}}x_{ij_{2}})\cdots
\\&&
(\sum\limits_{i=1}^{n}\overline {x}_{ij_{q}}x_{in_{p}})\cdots(\sum\limits_{i=1}^{n}\overline {x}_{in_{r}}x_{ik_{2}})\cdots(\sum\limits_{i=1}^{n}\overline {x}_{ik_{l}}x_{in_{r}}),
\end{array}
\label{0.3}
\end{equation}
where for sake of convenience, $x(t)$ is briefly denoted by $x$.
If we expand (\ref{0.3}), every expanded term can be expressed as
\[
\overline {x}_{i_{n_{1}}n_{1}}x_{i_{n_{1}}i_{2}}\cdots\overline {x}_{i_{i_{s}}i_{s}}x_{i_{i_{s}}n_{1}}\cdots\overline {x}_{i_{n_{p}}n_{p}}x_{i_{n_{p}}j_{2}}\cdots\overline {x}_{i_{j_{q}}j_{q}}x_{i_{j_{q}}n_{p}}\cdots\overline {x}_{i_{n_{r}}n_{r}}x_{i_{n_{r}}k_{2}}\cdots\overline {x}_{i_{k_{l}}k_{l}}x_{i_{k_{l}}n_{r}}\cdots,
\]
where $\overline {x}_{i_{n_{1}}n_{1}}x_{i_{n_{1}}i_{2}}$ is the $i_{n_{1}}$ form of $\sum\limits_{i=1}^{n}\overline {x}_{in_{1}}x_{ii_{2}}$ and similar for the rest.
If $i_{n_{1}}=i_{i_{2}}$, the product
\begin{equation}
\overline {x}_{i_{n_{1}}n_{1}}x_{i_{n_{1}}i_{2}}\overline {x}_{i_{i_{2}}i_{2}}x_{i_{i_{2}}i_{3}}\cdots=(x_{i_{n_{1}}i_{2}}\overline {x}_{i_{i_{2}}i_{2}})\overline {x}_{i_{n_{1}}n_{1}}x_{i_{i_{2}}i_{3}}\cdots,
\label{0.4}
\end{equation}
 is one of terms of $<\sigma>$. Thus,
we have
\[
\sigma_{1}^{\ast}=(n_{1}i_{3}\cdots i_{s}), \sigma_{1}^{\nabla}=(i_{2}),
\]
and
\[
\sigma^{\nabla}=\sigma_{1}^{\ast}\cdots \sigma_{1}^{\nabla} \cdots.
\]
We see that
\begin{equation}
\overline {x}_{i_{n_{1}}n_{1}}x_{i_{i_{2}}i_{3}}\cdots x_{i_{n_{1}}i_{2}}\overline {x}_{i_{i_{2}}i_{2}}\cdots
\label{0.5}
\end{equation}
is a term of $<\sigma^{\nabla}>$.
It is obvious that (\ref{0.4}) is equivalent to (\ref{0.5})
and $\varepsilon(\sigma)=-\varepsilon(\sigma^{\nabla})$. Thus, the term (\ref{0.4}) and the term (\ref{0.5}) canceled each other out
%$<\sigma>+<\sigma^{\nabla}>=0$,
when $i_{n_{1}}=i_{i_{2}}$. After the terms with $i_{n_{1}}=i_{i_{2}}$ canceled out, in the rest terms, $i_{n_{1}}$ is different from $i_{i_{2}}$.

If $i_{n_{1}}=i_{i_{\omega}}$, $i_{\omega}\in\{i_{3},\cdots,i_{s}\}$. Corresponding to the product
\begin{equation}
\overline {x}_{i_{n_{1}}n_{1}}x_{i_{n_{1}}i_{2}}\overline {x}_{i_{i_{2}}i_{2}}x_{i_{i_{2}}i_{3}}\cdots\overline {x}_{i_{i_{\omega-1}}i_{\omega-1}}x_{i_{i_{\omega-1}}i_{\omega}}\overline {x}_{i_{i_{\omega}}i_{\omega}}x_{i_{i_{\omega}}i_{\omega+1}}\overline {x}_{i_{i_{\omega+1}}i_{\omega+1}}x_{i_{i_{\omega+1}}i_{\omega+2}}\cdots,
\label{0.6}
\end{equation}
there is term expressed by
\begin{equation}
\overline {x}_{i_{n_{1}}n_{1}}x_{i_{i_{\omega}}i_{\omega}}\overline x_{i_{i_{\omega-1}}i_{\omega}}{x}_{i_{i_{\omega-1}}i_{\omega-1}}\cdots\overline x_{i_{i_{2}}i_{3}}
 {x}_{i_{i_{2}}i_{2}}\overline {x}_{i_{n_{1}}i_{2}}x_{i_{i_{\omega}}i_{\omega+1}}\overline {x}_{i_{i_{\omega+1}}i_{\omega+1}}x_{i_{i_{\omega+1}}i_{\omega+2}}\cdots,
 \label{0.7}
\end{equation}
which is one term of $<\sigma^{\diamond}>$.
\[
\sigma^{\diamond}=(n_{1},i_{\omega},i_{\omega-1},\cdots,i_{3},i_{2},i_{\omega+1},\cdots,i_{s})\cdots(n_{p},j_{2},\cdots,j_{q})\cdots(n_{r},k_{2},\cdots,k_{t})
=\sigma_{1}^{\diamond}\cdots\sigma_{p}\cdots\sigma_{r}.
\]
Note that
\[
x_{i_{n_{1}}i_{2}}\overline {x}_{i_{i_{2}}i_{2}}x_{i_{i_{2}}i_{3}}\cdots\overline {x}_{i_{i_{\omega-1}}i_{\omega-1}}x_{i_{i_{\omega-1}}i_{\omega}}\overline {x}_{i_{i_{\omega}}i_{\omega}}=\overline{x_{i_{i_{\omega}}i_{\omega}}\overline x_{i_{i_{\omega-1}}i_{\omega}}{x}_{i_{i_{\omega-1}}i_{\omega-1}}\cdots\overline x_{i_{i_{2}}i_{3}}
 {x}_{i_{i_{2}}i_{2}}\overline {x}_{i_{n_{1}}i_{2}}}.
\]
Thus, the product (\ref{0.6}) plus the product (\ref{0.7}) (for the sake of convenience, denoted by $ (\ref{0.6})+ (\ref{0.7})$) leads to
\[
 (\ref{0.6})+ (\ref{0.7})=2\Re(x_{i_{n_{1}}i_{2}}\overline {x}_{i_{i_{2}}i_{2}}x_{i_{i_{2}}i_{3}}\cdots\overline {x}_{i_{i_{\omega-1}}i_{\omega-1}}x_{i_{i_{\omega-1}}i_{\omega}}\overline {x}_{i_{i_{\omega}}i_{\omega}})\overline {x}_{i_{n_{1}}n_{1}}x_{i_{i_{\omega}}i_{\omega+1}}\cdots.
\]
Let $i_{m}$ be the biggest number in a set $\{i_{2},i_{3},\cdots,i_{\omega-1},i_{\omega}\}$, then we have
\[
\sigma_{1}^{\ast}=(n_{1}i_{\omega+1}\cdots i_{s}), \sigma_{1}^{\nabla}=(i_{m}i_{m+1}\cdots i_{\omega}i_{2}i_{3}\cdots i_{m-1}),\sigma_{1}^{\triangle}=(i_{m}i_{m-1}\cdots i_{3}i_{2}i_{\omega}\cdots i_{m+1}),
\]
and
\[
\sigma^{\nabla}=\sigma_{1}^{\ast}\cdots \sigma_{1}^{\nabla}\cdots , \sigma^{\triangle}=\sigma_{1}^{\ast}\cdots\sigma_{1}^{\triangle}\cdots.
\]
We see that
\begin{equation}
\overline {x}_{i_{n_{1}}n_{1}}x_{i_{i_{\omega}}i_{\omega+1}}\overline {x}_{i_{i_{\omega+1}}i_{\omega+1}}x_{i_{i_{\omega+1}}i_{\omega+2}}\cdots \overline x_{i_{i_{m}}i_{m}}{x}_{i_{i_{m}}i_{m+1}}\cdots\overline x_{i_{i_{\omega}}i_{\omega}}{x}_{i_{i_{\omega}}i_{2}}\cdots\overline x_{i_{i_{m-1}}i_{m-1}}{x}_{i_{i_{m-1}}i_{m}}\cdots
\label{0.8}
\end{equation}
is one term of $<\sigma^{\nabla}>$ and
\begin{equation}
\overline {x}_{i_{n_{1}}n_{1}}x_{i_{i_{\omega}}i_{\omega+1}}\overline {x}_{i_{i_{\omega+1}}i_{\omega+1}}x_{i_{i_{\omega+1}}i_{\omega+2}}\cdots \overline x_{i_{i_{m}}i_{m}}{x}_{i_{i_{m}}i_{m-1}}\cdots\overline x_{i_{i_{2}}i_{2}}{x}_{i_{i_{2}}i_{\omega}}\cdots\overline x_{i_{i_{m+1}}i_{m+1}}{x}_{i_{i_{m+1}}i_{m}}\cdots
\label{0.9}
\end{equation}
is one term of $<\sigma^{\triangle}>$. In view of $\sigma_{1}^{\nabla}=\overline{\sigma_{1}^{\triangle}}$, we have
\[
(\ref{0.8})+(\ref{0.9})=2\Re(\overline x_{i_{i_{m}}i_{m}}{x}_{i_{i_{m}}i_{m+1}}\cdots\overline x_{i_{i_{\omega}}i_{\omega}}{x}_{i_{i_{\omega}}i_{2}}\cdots\overline x_{i_{i_{m-1}}i_{m-1}}{x}_{i_{i_{m-1}}i_{m}})\overline {x}_{i_{n_{1}}n_{1}}x_{i_{i_{\omega}}i_{\omega+1}}\cdots.
\]
By $\Re(ab)=\Re(ba)$, $a,b$ are quaternions. Therefore
\[
\begin{array}{lll}
&&
\Re(x_{i_{n_{1}}i_{2}}\overline {x}_{i_{i_{2}}i_{2}}x_{i_{i_{2}}i_{3}}\cdots\overline {x}_{i_{i_{\omega-1}}i_{\omega-1}}x_{i_{i_{\omega-1}}i_{\omega}}\overline {x}_{i_{i_{\omega}}i_{\omega}})
\\&
=
&\Re(\overline x_{i_{i_{m}}i_{m}}{x}_{i_{i_{m}}i_{m+1}}\cdots\overline x_{i_{i_{\omega}}i_{\omega}}{x}_{i_{i_{\omega}}i_{2}}\cdots\overline x_{i_{i_{m-1}}i_{m-1}}{x}_{i_{i_{m-1}}i_{m}}).
\end{array}
\]
It is easy to see that $\varepsilon(\sigma)=\varepsilon(\sigma^{\diamond})=-\varepsilon(\sigma^{\nabla})=-\varepsilon(\sigma^{\triangle})$. Therefore, the terms $(\ref{0.6})+(\ref{0.7})$ and $(\ref{0.8})+(\ref{0.9})$ canceled each other out when $i_{n_{1}}=i_{i_{\omega}}$, $i_{\omega}\in\{i_{3},\cdots,i_{s}\}$.  After such terms with $i_{n_{1}}=i_{i_{\omega}}$ canceled out, in the rest terms, $i_{n_{1}}$ is different from $i_{i_{2}},i_{i_{3}},\cdots,i_{i_{s}},\cdots,i_{n_{p}},\cdots,i_{j_{q}},\\\cdots,i_{n_{r}},\cdots,i_{k_{l}}$.
%By similar analysis, we know that $i_{i_{\omega}}$ is different from $i_{i_{\omega+1}},i_{i_{\omega+2}}$ $\cdots,i_{i_{s}},\cdots,i_{n_{p}}\cdots,i_{j_{q}}\cdots i_{n_{r}}\cdots,i_{k_{l}}$.
The rest can be done in the same manner, $i_{j_{\omega}}$ is different from $i_{j_{\omega+1}},i_{j_{\omega+2}},\cdots,\\i_{j_{q}},\cdots,i_{n_{p}},\cdots,i_{i_{j_{q}}},\cdots ,i_{n_{r}},\cdots,i_{k_{l}}$. Finally we can prove that $i_{n_{1}},i_{i_{2}},i_{i_{3}},\cdots,i_{i_{s}},\cdots,i_{n_{p}},\cdots,\\i_{i_{j_{q}}},\cdots i_{n_{r}},\cdots,i_{k_{l}}$ are distinct from each other in the rest terms.
According to $\Re{a}_{ii}\in \mathbb{R}$. Obviously, $i_{n_{1}},i_{i_{2}},i_{i_{3}},\cdots,i_{i_{s}},\cdots,i_{n_{p}},\cdots,i_{i_{j_{q}}},\cdots ,i_{n_{r}},\cdots,i_{k_{l}}$ are distinct from each other in the rest terms of $A_{1}(t),A_{2}(t),\cdots,A_{n}(t)$.

For $W_{QDE}^P(t)$, letting
\[
\overline {x}_{i_{n_{1}}n_{1}}x_{i_{n_{1}}i_{2}}\cdots\overline {x}_{i_{i_{s}}i_{s}}x_{i_{i_{s}}n_{1}}\cdots\overline {x}_{i_{n_{p}}n_{p}}x_{i_{n_{p}}j_{2}}\cdots\overline {x}_{i_{j_{q}}n_{p}}x_{i_{j_{q}}n_{p}}\cdots\overline {x}_{i_{n_{r}}n_{r}}x_{i_{n_{r}}k_{2}}\cdots\overline {x}_{i_{k_{l}}k_{l}}x_{i_{k_{l}}n_{r}}=a
\]
be one term of the rest in $<\sigma>$. And for the rest terms of $A_{1}(t),A_{2}(t),\cdots,A_{n}(t)$, the coefficient of $a$ must be distinct from each other. Therefore, by the arbitrariness of $a$ and $<\sigma>$, we easily get
\begin{equation}
\sum\limits_{i}^{n}A_{i}(t)=\sum\limits_{i}^{n}(2\Re{a}_{ii}a+\cdots)=\sum\limits_{i}^{n}2\Re{a}_{ii}(a+\cdots)=\sum\limits_{i}^{n}2\Re{a}_{ii}W_{QDE}^P(t).
\label{0.10}
\end{equation}

(ii) For arbitrary $k,j$ ($k\neq j$), we will prove $\sum\limits_{m=1}^{n}B_{kj}^{m}(t)=0$. Without loss of generality, let $ k=1$,$j=2$.
For every permutation $\sigma\in S_{n}$
\[
\sigma=(n_{1}i_{2}\cdots i_{s})\cdots(n_{p}j_{2}\cdots j_{q})\cdots(n_{r} k_{2} \cdots k_{t})=\sigma_{1}\cdots\sigma_{p}\cdots\sigma_{r}.
\]
Consider $B_{12}^{j_{\omega}}(t)$ and $\sigma_{p}=(n_{p}j_{2}\cdots j_{\omega-1} j_{\omega} j_{\omega+1}\cdots j_{q})$. Every expanded terms of $<\sigma>$ can be expressed as
\begin{equation}
\overline {x}_{i_{n_{1}}n_{1}}x_{i_{n_{1}}i_{2}}\cdots\cdots\overline {x}_{i_{j_{\omega-1}}j_{\omega-1}}{x}_{i_{j_{\omega-1}}1}\overline x_{2j_{\omega}}\overline a_{12} x_{1j_{\omega+1}}
\overline {x}_{i_{j_{\omega+1}}j_{\omega+1}}{x}_{i_{j_{\omega+1}}j_{\omega+2}}\cdots\cdots\overline {x}_{i_{k_{l}}k_{l}}{x}_{i_{k_{l}}n_{r}},
\label{0.11}
\end{equation}
or
\begin{equation}
\overline {x}_{i_{n_{1}}n_{1}}x_{i_{n_{1}}i_{2}}\cdots\cdots\overline {x}_{i_{j_{\omega-1}}j_{\omega-1}}{x}_{i_{j_{\omega-1}}1}\overline x_{1j_{\omega}}a_{12} x_{2j_{\omega+1}}
\overline {x}_{i_{j_{\omega+1}}j_{\omega+1}}{x}_{i_{j_{\omega+1}}j_{\omega+2}}\cdots\cdots\overline {x}_{i_{k_{l}}k_{l}}{x}_{i_{k_{l}}n_{r}}.
\label{0.12}
\end{equation}
From the proof of (i), we easily obtain that $i_{n_{1}},\cdots,i_{j_{\omega-1}},i_{j_{\omega+2}},\cdots,i_{k_{l}}$ are distinct from each other in the rest terms.
Because $i_{n_{1}},\cdots,i_{j_{\omega-1}},i_{j_{\omega+2}},\cdots,i_{k_{l}}$ are distinct from each other in the rest forms, there is $i_{z}=1$ or $i_{z}=2$ or $i_{z_{1}}=1,i_{z_{2}}=2$, $i_{z},i_{z_{1}},i_{z_{2}}\in\{i_{n_{1}},\cdots,i_{j_{\omega-1}},i_{j_{\omega+2}},\cdots,i_{n_{1}}\}$.

For $i_{z}=1$ or $i_{z}=2$ or $i_{z_{1}}=1,i_{z_{2}}=2$, since (\ref{0.11}) contains $\overline x_{1j_{\omega}}$ and  $x_{2j_{\omega+1}}$,   (\ref{0.11}) can be canceled out from the proof of (i). And so is (\ref{0.12}). Then, we obtain $\sum\limits_{m=1}^{n}B_{12}^{m}(t)=0$, and $\sum\limits_{m=1}^{n}B_{kj}^{m}(t)=0$. Therefore, according to (\ref{0.2}) and (\ref{0.10}),  we obtain
\[
\frac{d}{dt}W_{QDE}^P(t)=\sum\limits_{m=1}^{n}A_{m}(t)=\sum\limits_{i}^{n}2\Re{a}_{ii}(t)W_{QDE}^P(t)
=[trA(t)+trA^{+}(t)]W_{QDE}^P(t).
\]
Integration above equation over $[t_0,t]$ follows the Liouville formula.

We need a lemma from (Theorem 2.10 \cite{Ky-AMC}).

 \begin{lemma}\label{le4.2}
A quaternionic matrix $M$ is invertible if and only if $\rm{ddet}_P M \neq0$.
  \end{lemma}

\begin{proposition}\label{Pr4.1}
If $W_{QDE}^P(t)=0$ at some $t_0$ in $I$ then  $x_1(t),x_2(t),\cdots,x_n(t)$ are right dependent on $I$.
\end{proposition}

 {\bf Proof.} From Liouville formula, we have
 \[
 W_{QDE}^P(t_0)=0,\,\,\,\,\, \mbox{implies} \,\,\,\,\, W_{QDE}^P(t)=0,\,\,\,\, \mbox{for any}\,\,\,\,t\in I.
 \]
According to Lemma \ref{le4.2}, the quaternionic matrix $M(t)$ is not invertible on $I$. Hence, the linear system
 \[
 M(t) r=0,\,\,\,\, or \,\,\,\, (x_1(t),x_2(t),\cdots,x_n(t))r=0,     \,\,\,\, r=(r_1,r_2,\cdots,r_n)^{T}\in \mathbb{H}^{n},
 \]
 has a non-zero solution. Consequently, the $n$ solution $x_1(t),x_2(t),\cdots,x_n(t)$  are right dependent on $I$.
 %Similarly, it is easy to prove the left linear dependence.

 From Lemma \ref{Th4.2}, \ref{Th4.3} and Theorem \ref{Th4.4}, we immediately have

 \smallskip

\begin{lemma}\label{Th4.5}
Let $A(t)$ in Eq.$($\ref{2.1}$)$ be continuous functions defined on an interval $t\in I$. $n$ solutions $x_1(t),x_2(t),\cdots,x_n(t)$ of Eq.$($\ref{2.1}$)$ on $I$ are right dependent on $I$ if and only if the Wronskian $W_{QDE}^P(t_0)=0$.
\end{lemma}

\smallskip
 Now, we present two important results on the structure of the general solution.
\begin{lemma}\label{Th4.6}
There are $n$ independent solutions $x_1(t), x_2(t), \cdots , x_n(t) $ of Eq.$($\ref{2.1}$)$ associated with the initial value problem \eqref{2.2}.
\end{lemma}

 The proof is similar to that of ODEs. And it is not difficult to prove the following theorem on the general solutions by above lemmas.

\begin{theorem}\label{Th4.7}
$($Structure of the general solution$)$ If $x_1(t), x_2(t), \cdots , x_n(t) $ are n independent solutions of Eq.$($\ref{2.1}$)$, then each solution of Eq.$($\ref{2.1}$)$ is expressed as
 \begin{equation}
 x(t)=x_1(t) r_1 + x_2(t) r_2 + \cdots + x_n(t) r_n,
 \label{sgs}
 \end{equation}
where $r_{1}, r_{2},\cdots, r_{n}$ are undetermined quaternionic constants. The set of all the solutions consists of {\bf a free right-module}.
 \end{theorem}

\section{Fundamental Matrix and Solution to QDEs}

%\begin{proposition}\label{Pr5.4}
%If $A(t)$ is differentiable and if $A(t)A'(t) = A'(t)A(t)$, it
%follows from that
%\begin{equation}
%[\exp A(t) ]' = [\exp A(t) ] A'(t).
%\label{EX2}
%\end{equation}
%\end{proposition}

\begin{definition}\label{De5.1}
A {\em solution matrix} of Eq.$($\ref{2.1}$)$ was denoted by
\[
M(t)=(x_1(t),x_2(t),\cdots,x_n(t))
=
\left ( \begin{array}{ll}
  {x} _{11}(t)  & {x} _{12}(t) \cdots {x} _{1n}(t)
 \\
 {x} _{21}(t)  & {x} _{22}(t) \cdots {x} _{2n}(t)
 \\
 & \cdots
 \\
  {x} _{n1}(t)  & {x} _{n2}(t) \cdots {x} _{nn}(t)
\end{array}
\right ),
\]
where $x_1(t),x_2(t),\cdots,x_n(t)$ be any $n$ solutions of Eq.$($\ref{2.1}$)$ on $I$. Moreover, if $x_1(t),x_2(t),\cdots,x_n(t)$ are independent, we call it a {\em fundamental matrix} of Eq.$($\ref{2.1}$)$. In particular, if $M(t_0)$ is an identity, then we call it a {\em normal fundamental matrix}.
\end{definition}

%From Theorem \ref{Th4.3}, we know that if $M(t)$ is a fundamental matrix of Eq.(\ref{2.1}), the Wronskian determinant $W_{QDE}^P(t)\neq0$. By Theorem \ref{Th4.7}, we have

\begin{theorem}\label{Th5.1}
A general solution $x(t)$ of  Eq.$($\ref{2.1}$)$ can be rewritten as
\[
x(t)=M(t) q,
\]
where $M(t)$ is a fundamental matrix, $q\in   \mathbb{H}^{n} $ is a constant quaternionic vector. In particular, if $A(t)\equiv A$ (a constant quaternion matrix), then
\[
x(t)=\exp\{A t\} q.
\]
For any given initial value $x(t_0)=x^0$, the corresponding solution to the initial value is
\[
x(t)=M(t)M^{-1}(t_0)x^0.
\]
In particular, if $A(t)\equiv A$ (a constant quaternion matrix), then
\[
x(t)=\exp\{A( t-t_{0})\} x^0.
\]
\end{theorem}

\begin{remark}\label{Th5.2}
  A solution matrix $M(t)$ of Eq.$($\ref{2.1}$)$ on $I$ is a fundamental matrix if and only if $\rm{ddet}_P M(t_0)=\neq 0$ $($or $W_{QDE}^P(t_0)\neq 0$$)$ on $I$.
  \end{remark}

\smallskip

In particular, now we consider $A(t)\equiv A$ (a constant quaternion matrix), that is, QDEs with constant quaternionic matrix
\begin{equation}
\dot{x}(t) = A x(t).
\label{Constant}
\end{equation}

\begin{theorem}\label{Th5.4}
Assume that the commutivity property
\begin{equation}
a_i(t)   \int_{t_0}^{t} a_i(s) ds\  =  \int_{t_0}^{t} a_i(s) ds  a_i(t)
\label{Condition1}
\end{equation}
holds. Then the fundamental matrix of the diagonal homogenous system
\begin{equation}
\left ( \begin{array}{ll}
 \dot{x} _{1}(t)
 \\
\dot{x} _{2}(t)
\\
\,\,\,\vdots
\\
\dot{x} _{n}(t)
\end{array}
\right )
=
\left ( \begin{array}{ccc}
{a} _{1}(t) & \,\,\,\,0\,\,\,\,\,\,\,\cdots \,\,\,\, 0
\\
0 & {a} _{2}(t) \,\cdots \,\,\,\,0
\\
\vdots & \,\,\vdots \,\,\,\,\,\, \,\ddots\,\,\,\,\,\, \vdots
\\
0 & \,\,\,\,\,\,\,\,0 \,\,\,\,\,\,\cdots  \,\,{a} _{n}(t)
\end{array}
\right )
\left ( \begin{array}{ll}
{x} _{1}(t)
\\
{x} _{2}(t)
\\
\,\,\,\vdots
\\
{x} _{n}(t)
\end{array}
\right ).
\label{diag}
\end{equation}
can be chosen as
\[
M(t)=
\left ( \begin{array}{ccc}
\exp\{ \int_{t_0}^{t} a_1(s) ds\}   & \,\,\,\,\,\,\,\,\,\,\,\,\,\,\,\,\,\,\,\,0 \,\,\,\,\,\,\,\,\,\,\,\,\,\,\,\,\,\,\,\,\cdots \,\,\,\,\,\,\,\,\,\,\,\,\,\,\,\,\,\,\,\,\,\,\,\, 0
 \\
0 &\exp\{ \int_{t_0}^{t} a_2(s) ds\} \,\cdots \,\,\,\,\,\,\,\,\,\,\,\,\,\,\,\,\,\,\,\,\,\,\,\,0
\\
\vdots & \,\,\,\,\,\,\,\,\,\,\,\,\,\,\,\,\,\,\,\,\vdots \,\,\,\,\,\,\,\,\,\,\,\,\,\,\,\,\,\,\,\, \ddots\,\,\,\,\,\, \,\,\,\,\,\,\,\,\,\,\,\,\,\,\,\,\,\,\vdots
\\
0 & \,\,\,\,\,\,\,\,\,\,\,\,\,\,\,\,\,\,\,\,\,\,\,\,\,\,\,\,\,\,\,\,\,\,\,\,\,\,\,0 \,\,\,\,\,\,\,\,\,\,\,\,\,\,\,\,\,\,\,\cdots \,\,\,\,\exp\{ \int_{t_0}^{t} a_n(s) ds\}
\end{array}
\right ).
\]
 \end{theorem}

\section{\bf Algorithm for computing fundamental matrix %of linear QDES with constant coefficients
}

%Since the fundamental matrix plays great role in solving QDEs, in
Two algorithms for computing fundamental matrix of linear QDEs with constant coefficients will given in this section.

\subsection{ Method 1: using expansion of $\exp \{At\}$}

%It is easy to see that $\exp\{A t\}$ is a fundamental matrix of Eq.(\ref{Constant}).
%So if the coefficient matrix is not very complicate, we can use the definition of $\exp\{A t\}$ to compute fundamental matrix of linear QDEs with constant coefficients.

\begin{theorem}\label{Th6.1} If $A=diag( \lambda_1,\lambda_2,\cdots,\lambda_n)\in \mathbb{H}^{n\times n}$ is a diagonal matrix, then
\[
\exp \{At \}= \left ( \begin{array}{ccc}
\exp\{\lambda_1 t\} & \,\,\,\,\,\,\,\,0 \,\,\,\,\,\,\,\,\,\,\,\,\cdots \,\,\,\,\,\,\,\, 0
\\
0 & \exp\{\lambda_2 t\} \cdots \,\,\,\,\,\,\,\,0
\\
\vdots &\,\,\,\,\,\,\,\,\, \vdots \,\,\,\,\,\,\,\,\,\,\, \ddots\,\,\,\,\,\,\,\,\,\vdots
\\
0 & \,\,\,\,\,\,\,\,\,\,\,\,\,\,\,\,\,\,\,\,0\, \,\,\,\,\,\,\,\,\,\,\cdots \exp\{ \lambda_n t\}
\end{array}
\right ).
\]
\end{theorem}

{\bf Proof.} By the expansion,
\[
\begin{array}{cccc}
\exp \{ At \}
&=&
E+ \left ( \begin{array}{ccc}
\lambda_1 & \,0\,\,\,\cdots \,\,\,\, 0
\\
0 & \lambda_2 \,\cdots \,\,\,\,0
\\
\vdots & \vdots \,\,\,\,\ddots\,\,\,\, \vdots
\\
0 & \,\,\,0 \,\,\,\cdots \,\,\,\lambda_n

\end{array}
\right )\frac{t}{1!}
+
 \left ( \begin{array}{ccc}
\lambda_1 & \,0\,\,\,\cdots \,\,\,\, 0
\\
0 & \lambda_2 \,\cdots \,\,\,\,0
\\
\vdots & \vdots \,\,\,\,\ddots\,\,\,\, \vdots
\\
0 & \,\,\,0 \,\,\,\cdots \,\,\,\lambda_n

\end{array}
\right )^2\frac{t^2}{2!}
+\cdots
\\
&=&\left(\begin{array}{ccc}
\exp\{\lambda_1 t\} & \,\,\,\,\,\,\,\,0 \,\,\,\,\,\,\,\,\,\,\,\,\cdots \,\,\,\,\,\,\,\, 0
\\
0 & \exp\{\lambda_2 t\} \cdots \,\,\,\,\,\,\,\,0
\\
\vdots &\,\,\,\,\,\,\,\,\, \vdots \,\,\,\,\,\,\,\,\,\,\, \ddots\,\,\,\,\,\,\,\,\,\vdots
\\
0 & \,\,\,\,\,\,\,\,\,\,\,\,\,\,\,\,\,\,\,\,0\, \,\,\,\,\,\,\,\,\,\,\cdots \exp\{ \lambda_n t\}
\end{array}
\right )
.
\end{array}
\]

If we can divide the matrix to some simple ones and use the expansion to compute the fundamental matrix.
\[
A=diag A+ N,
\]
where $N$ is a nilpotent matrix. That is, $N^n=0$ and $n$ is a finite number.

\smallskip

\noindent {\bf Example 5.1} Find a fundamental matrix of the following QDES
\[
\dot {x}= Ax=
\left ( \begin{array}{ccccc}
\lambda & 1 &0 \cdots 0 & 0
\\
0 & \lambda &1 \cdots 0 & 0
\\
\,\,& \,\,& \,\,\,\,\,\,\, \ddots \,\,\ddots  &
\\
0 & 0 &0\cdots \lambda & 1
\\
0 & 0 &0\cdots 0& \lambda
\end{array}
\right )x,\,\,\,\, x=(x_1,x_2,\cdots,x_k)^T.
\]

{\bf Answer.} We see that
$A=\lambda E +B$. Noticing that $(\lambda E) B=B(\lambda E)$, by Theorem \ref{Th6.1}, we have $\exp\{At\}=\exp\{\lambda E t\}\cdot\exp\{Bt\}$, where
\[
B=\left ( \begin{array}{ccccc}
0 & 1 &0 \cdots 0 & 0
\\
0 & 0 &1 \cdots 0 & 0
\\
\,\,& \,\,& \,\,\,\,\,\,\, \ddots \,\,\ddots  &
\\
0 & 0 &0 \cdots 0 & 1
\\
0 & 0 &0 \cdots 0 & 0
\end{array}
\right ).
\]
Note that $B$ is a nilpotent matrix. That is, $B^k=0$, we get
\[
\exp\{Bt\}=
\left ( \begin{array}{ccccc}
1 & t &\frac{t^{2}}{2!} \cdots \frac{t^{k-2}}{(k-2)!} & \frac{t^{k-1}}{(k-1)!}
\\
0 & 1 &t \cdots \frac{t^{k-3}}{(k-3)!} & \frac{t^{k-2}}{(k-2)!}
\\
\,\,& \,\,& \,\,\,\,\, \ddots \,\,\ddots  &
\\
0 & 0 &0 \cdots 1 & t
\\
0 & 0 &0 \cdots 0 & 1
\end{array}
\right ).
\]
Then the fundamental matrix
\[
\exp\{At\}=\exp\{\lambda E t\}\cdot\exp\{Bt\}
=
\left ( \begin{array}{ccccc}
1 & t &\frac{t^{2}}{2!} \cdots \frac{t^{k-2}}{(k-2)!} & \frac{t^{k-1}}{(k-1)!}
\\
0 & 1 &t \cdots \frac{t^{k-3}}{(k-3)!} & \frac{t^{k-2}}{(k-2)!}
\\
\,\,& \,\,& \,\,\,\,\, \ddots \,\,\ddots  &
\\
0 & 0 &0 \cdots 1 & t
\\
0 & 0 &0 \cdots 0 & 1
\end{array}
\right )\exp\{ \lambda t\}.
\]

\subsection {Method 2:  eigenvalue and eigenvector theory}

The eigenvalues of quaternion matrices should be treated as left eigenvalues and right eigenvalues. Usually, they are different and not equal. They have no relations. The numbers of the eigenvalues are possible to be infinite. Thus, an eigenvalue $\theta$ is similar to $\lambda$ if
\[
\theta=\alpha^{-1} \lambda \alpha
\]
 \begin{remark}\label{Re6.1.1}
If $\theta$, $\lambda$ are two characteristic roots of $A$ and $\theta$ is similar to $\lambda$, for any the eigenvector $q$ of $\theta$, there exists an eigenvector $q'$ of $\lambda$ such that $q,q'$ are dependent.
\end{remark}

From the definition of fundamental matrix, we have

\begin{theorem}\label{Th6.2}
If the matrix $A$ has $n$ independent eigenvectors $q_1,q_2,\cdots,q_n$, corresponding to the eigenvalues $\lambda_1,\lambda_2,\cdots,\lambda_n $ $($$\lambda_i$ and $\lambda_j$ can be similar$)$, then the fundamental matrix of Eq.$($\ref{Constant}$)$ can be chosen as
\[
M(t)=(q_1e^{\lambda_1 t},q_2e^{\lambda_2 t},\cdots,q_n e^{\lambda_n t}).
\]

\end{theorem}

The proof is similar to Theorem 6.5 in \cite{QDE1}.

%{\bf Proof.} From above discussion, we know $q_1e^{\lambda_1 t},q_2e^{\lambda_2 t},\cdots,q_n e^{\lambda_n t}$ are $n$ solution of Eq.(\ref{Constant}). Thus, $M(t)$ is a solution matrix of Eq.(\ref{Constant}). Moreover, by using the independence of $q_1,q_2,\cdots,q_n$, we have
%\[
%{\rm ddet}M(0)={\rm ddet}(q_1,q_2,\cdots,q_n) \neq 0.
%\]
%Therefore, $M(t)$ is a fundamental matrix of Eq.(\ref{Constant}).

Then by Lemma 6.6 in \cite{QDE1}, we have the corollary

\begin{corollary}\label{Co6.1} If the matrix $A$ has $n$ distinct eigenvalues $\lambda_1,\lambda_2,\cdots,\lambda_n$, no two of which are
similar, then the fundamental matrix of Eq.$($\ref{Constant}$)$ can be chosen as
\[
M(t)=(q_1e^{\lambda_1 t},q_2e^{\lambda_2 t},\cdots,q_n e^{\lambda_n t}).
\]
\end{corollary}

%Some results from \cite{Bren} (Theorem 11, Theorem 12).

%\begin{lemma}\label{Le6.2}
%If $A$ is in triangular form, then every diagonal element is a
%characteristic root.
%\end{lemma}

%\begin{lemma}\label{Le6.3}
%Let a matrix of quaternion be in triangular form. Then the only
%characteristic roots are the diagonal elements $($and the numbers similar to them$)$.
%\end{lemma}

%\begin{lemma}\label{Le6.4}
% Similar matrices have the same characteristic roots.
% \end{lemma}

%\begin{lemma}\label{Le6.5}
 %\cite{Bren} $($Theorem 2$)$ Every matrix of quaternion can be transformed into triangular form by a unitary matrix.
 %\end{lemma}

When the matrix $A$ has $n$ distinct (simple) eigenvalues, we have shown how to construct a fundamental matrix in \cite{QDE1}. But how to construct the fundamental matrix when the matrix $A$ has multiple eigenvalues? Next section is devoted to answering this question.

\section{System with multiple eigenvalues}

In this section, we will give an algorithm to construct fundamental matrix when system have the multiple eigenvalues. There are two cases. One case is that the numbers of eigenvectors are equal to the dimension of the system. The other case is that the numbers of eigenvectors less than the dimension of the system (that is to say, not enough eigenvectors). It should be noted that two similar eigenvalues can be seen as an eigenvalue with the multiplicity of two.

\subsection{Multiple eigenvalues with enough eigenvectors}

%In this subsection, we give two examples. Example 2 has two identity eigenvalues. Example 3 has two eigenvalues which is similar.

\noindent {\bf Example 6.1} Find a fundamental matrix of the following QDEs
\begin{equation}
\dot {x}= \left ( \begin{array}{ll}
 \boldsymbol{j}  &  \boldsymbol{i}
 \\
 0  &  \boldsymbol{j}
\end{array}
\right )x,\,\,\,\, x=(x_1,x_2)^T.
\label{Ex6.2}
\end{equation}

{\bf Answer:} From Lemma 6.8 and Lemma 6.9 in  \cite{QDE1}, we see that $\lambda_{1,2}=\boldsymbol{j}$.
To find the eigenvector of $\lambda_{1,2}=\boldsymbol{j}$, we consider the following equation
\[
Aq=q\lambda_{1,2},
\]
that is
\begin{equation}
 \left\{ \begin{array}{ccc}
 \boldsymbol{j} q_1 + \boldsymbol{i} q_2 & =& q_1 \boldsymbol{j},
 \\
\boldsymbol{j} q_2 &=& q_2 \boldsymbol{j}.
\end{array}
\right.
\label{EV1}
\end{equation}
From the second equation of (\ref{EV1}), if we take $q_2=0$. Substituting it into the first equation of (\ref{EV1}), we can take
$q_1=1$. So we obtain one eigenvector as
\[
\nu_1=\left ( \begin{array}{ll}
 q_1
 \\
q_2
\end{array}
\right )
=\left ( \begin{array}{ll}
 1
 \\
0
\end{array}
\right ).
\]
If we take $q_2=1$, substituting it into the first equation of (\ref{EV1}), we can take
$q_1=-\frac{\boldsymbol{k}}{2}$. So we get another eigenvector as
\[
\nu_2=\left ( \begin{array}{ll}
 q_1
 \\
q_2
\end{array}
\right )
=\left ( \begin{array}{ll}
-\frac{\boldsymbol{k}}{2}
 \\
\,\,\,1
\end{array}
\right ).
\]
Since
\[
\begin{array}{lll}
{\rm ddet}(\nu_1,\nu_2)
=
{\rm ddet}\left ( \begin{array}{ll}
1 & -\frac{\boldsymbol{k}}{2}
 \\
0 & \,\,\,1
\end{array}
\right )
=
{\rm det}\Big[\left ( \begin{array}{ll}
1 & 0
 \\
\frac{\boldsymbol{k}}{2} & 1
\end{array}
\right )\left ( \begin{array}{ll}
1 & -\frac{\boldsymbol{k}}{2}
 \\
0 &\,\,\, 1
\end{array}
\right )
\Big]=1
\neq0,
\end{array}
\]
the eigenvectors $\nu_1$ and $\nu_2$ are independent. Taking
 \[
 M(t)=(\nu_1 e^{\boldsymbol{j} t},\nu_2 e^{\boldsymbol{j} t})=\left ( \begin{array}{ll}
 e^{\boldsymbol{j} t}  & -\frac{\boldsymbol{k}}{2} e^{\boldsymbol{j} t}
 \\
0  &  \,\,\,\,\,e^{\boldsymbol{j} t}
\end{array}
\right ).
 \]
From Theorem \ref{Th6.2}, $M(t)$ is a fundamental matrix. In fact, by definition of fundamental matrix, we can also verify that $M(t)$ is a fundamental matrix of Eq.(\ref{Ex6.2}). To show this fact, firstly, we show that $M(t)$ is a solution matrix of Eq.(\ref{Ex6.2}).
 Let $\phi_1(t)=\nu_1 e^{\lambda_{1,2} t}$ and $\phi_2(t)=\nu_2 e^{\lambda_{1,2} t}$, then
\[
\dot{\phi}_1(t)
=\left ( \begin{array}{ll}
\boldsymbol{j} e^{\boldsymbol{j} t}
 \\
\,\,0
\end{array}
\right )
=
\left ( \begin{array}{ll}
 \boldsymbol{j}  &  \boldsymbol{i}
 \\
 0  &  \boldsymbol{j}
\end{array}
\right )
\left ( \begin{array}{ll}
e^{\boldsymbol{j} t}
 \\
\,\,0
\end{array}
\right )
=
\left ( \begin{array}{ll}
 \boldsymbol{j}  &  \boldsymbol{i}
 \\
 0  &  \boldsymbol{j}
\end{array}
\right ) \phi_1(t),
\]
which implies that $\phi_1(t)$ is a solution of Eq.(\ref{Ex6.2}).
Similarly, we have
\[
\dot{\phi}_2(t)=\left ( \begin{array}{ll}
-\frac{\boldsymbol{k}}{2} \boldsymbol{j} e^{\boldsymbol{j} t}
 \\
\,\,\,\,\boldsymbol{j} e^{\boldsymbol{j} t}
\end{array}
\right )
=
\left ( \begin{array}{ll}
 \boldsymbol{j}  &  \boldsymbol{i}
 \\
 0  &  \boldsymbol{j}
\end{array}
\right )
\left ( \begin{array}{ll}
-\frac{\boldsymbol{k}}{2}
 \\
\,\,\,1
\end{array}
\right ) e^{\boldsymbol{j} t}
=
\left ( \begin{array}{ll}
 \boldsymbol{j}  &  \boldsymbol{i}
 \\
 0  &  \boldsymbol{j}
\end{array}
\right ) \phi_2(t),
\]
which implies that $\phi_2(t)$ is another solution of Eq.(\ref{Ex6.2}). Therefore,  $M(t)=(\phi_1(t),\phi_2(t))^T$ is a solution matrix of Eq.(\ref{Ex6.2}).

Secondly, in view of Theorem \ref{Th5.2} and the fact
\[
\begin{array}{lll}
\rm{ddet}_P M(t_0)=\rm{ddet}_P M(0)={\rm ddet}(\nu_1,\nu_2)
=
{\rm ddet}\left ( \begin{array}{ll}
1 & 1
 \\
0 & 1
\end{array}
\right )
=
{\rm det}\Big[\left ( \begin{array}{ll}
1 & 0
 \\
1 & 1
\end{array}
\right )\left ( \begin{array}{ll}
1 & 1
 \\
0 & 1
\end{array}
\right )
\Big]
\neq0,
\end{array}
\]
then $\rm{ddet}_P M(t)\neq0$. Therefore, $M(t)$ is a fundamental matrix of Eq. (\ref{Ex6.2}).

\subsection{Multiple eigenvalues with fewer eigenvectors}

For any $A\in\mathbb{H}^{n\times n}$, if we obtain double or multiple eigenvalues. This means that the number of independent eigenvectors might be less than the dimension of the system. So we may not get a fundamental matrix. We therefore have to discover how to find the ``missing solutions". In this case, first, we need to prove the following basic results.

Now we need a lemma from \cite{Zhangf2} (Theorem 5.4).

\begin{lemma}\label{Le6.6}
 Any $n\times n$ quaternion  matrix A has exactly $n$ $($right$)$ eigenvalues which are complex numbers with nonnegative imaginary parts.
 \end{lemma}

For $A\in\mathbb{H}^{n\times n}$, suppose that $\lambda_{1},\lambda_{2},\cdots,\lambda_{k}$ are distinct standard eigenvalues for $A$, the multiplicity of all the eigenvalues are $n_{1},n_{2},\cdots,n_{k}$ respectively, and $n_{1}+n_{2}+\cdots+n_{k}=n$. For the independent  eigenvectors $v_{1}^{j1},v_{1}^{j2},\cdots,v_{1}^{j r_{j}}$ associated with eigenvalue $\lambda_{j}$ of multiplicity $n_{j}$ ($r_{j}\leq n_{j}$). A set $\{v_{1}^{ji},v_{2}^{ji},\cdots,v_{m_{ji}}^{ji}\}$ $(m_{i}\leq n_{j})$ based on the eigenvector $v_{1}^{ji}$ ($A v_{1}^{ji}=\lambda_{j} v_{1}^{ji}$ and $i\in \{1,2,\cdots,r_{j}\}$) such that
\begin{equation}
\begin{array}{ccc}
 A v_{m_{ji}}^{ji}-v_{m_{ji}}^{ji}\lambda_{j}&=&v_{m_{ji}-1}^{ji},
\\
A v_{m_{ji}-1}^{ji}-v_{m_{ji}-1}^{ji}\lambda_{j}&=&v_{m_{ji}-2}^{ji},
\\
&\vdots &
\\
A v_{2}^{ji}-v_{2}^{ji}\lambda_{j}&=&v_{1}^{ji}.
\label{vl}
\end{array}
\end{equation}
Note that $N_{\lambda_{j}}^{m_{ji}}=\{\sum\limits^{m_{ji}}_{l=1} v_{l}^{ji}r_{l}|r_{l}\in\mathbb{H}\}$, and $N_{\lambda_{j}}^{m_{ji}}\subseteq\mathbb{H}^{n}$
is submodule.

Now we obtain the following result according to Lemma \ref{Lemma8}, Theorem 2 in \cite{Zhangf} and some basic theories of direct sum.

\begin{theorem}\label{Th6.4} For all submodules $N_{\lambda_{j}}^{m_{ji}}$, $i\in \{1,2,\cdots,r_{j}\}$ and $j\in \{1,2,\cdots,k\}$, then there exists the following decomposition
\begin{equation}
\mathbb{H}^{n}=N_{\lambda_{1}}^{m_{11}}\oplus N_{\lambda_{1}}^{m_{12}}\oplus\cdots\oplus N_{\lambda_{1}}^{m_{1 r_{1}}}\oplus N_{\lambda_{2}}^{m_{21}}\oplus\cdots\oplus N_{\lambda_{2}}^{m_{2 r_{2}}}
\oplus\cdots\oplus N_{\lambda_{k}}^{m_{k1}}\oplus\cdots\oplus N_{\lambda_{k}}^{m_{k r_{k}}},
\label{hn}
\end{equation}
where~$N_{\lambda_{j}}^{m_{ji}}=\{\sum\limits^{m_{ji}}_{l=1} v_{l}^{ji}r_{l}|r_{l}\in\mathbb{H}\}$.
\end{theorem}
From Theorem \ref{hn}, for any $u\in\mathbb{H}^{n}$, there exist the unique vectors $u_{1}^{m_{11}},\cdots,u_{1}^{m_{1 r_{1}}},u_{2}^{m_{21}},\\\cdots,u_{2}^{m_{2 r_{2}}},\cdots,$
$u_{k}^{m_{k1}},\cdots,u_{k}^{m_{k r_{k}}}$, where $u_{j}^{m_{ji}}\in N_{\lambda_{j}}^{m_{ji}}$, such that
\begin{equation}
u=u_{1}^{m_{11}}+\cdots+u_{1}^{m_{1 r_{1}}}+u_{2}^{m_{21}}+\cdots+u_{2}^{m_{2 r_{2}}}+\cdots+u_{k}^{m_{k1}}
+\cdots+u_{k}^{m_{k r_{k}}}=\sum\limits^{k}_{j=1}\sum\limits^{r_{j}}_{i=1}u_{j}^{m_{ji}}.
\label{uuj}
\end{equation}
By Eq.(\ref{uuj}) and $N_{\lambda_{j}}^{m_{ji}}$, any solution $x(t)=\exp\{A t\} \eta$ of Eq.(\ref{Constant}) can be represented by
\begin{equation}
\begin{array}{lll}
x(t)
 &=&
  (\exp\,At) \eta=(\exp\,At)\sum\limits^{k}_{j=1}\sum\limits^{r_{j}}_{i=1}u_{j}^{m_{ji}}=\sum\limits^{k}_{j=1}
  \sum\limits^{r_{j}}_{i=1}(\exp\,At)u_{j}^{m_{ji}}
%  \\
%  &=&
%  \sum\limits^{k}_{j=1}\sum\limits^{r_{j}}_{i=1}(\exp\,At)(v_{1}^{ji}r_{1}+v_{2}^{ji}r_{2}+\cdots+v_{m_{ji}}^{ji}r_{m_{ji}})
%  \\
 % &=&
%  \sum\limits^{k}_{j=1}\sum\limits^{r_{j}}_{i=1}\sum\limits^{+\infty}_{s=0}(A t)^{s}(v_{1}^{ji}r_{1}+v_{2}^{ji}r_{2}+\cdots+v_{m_{ji}}^{ji}r_{m_{ji}})
%  \\
%  &=&
%  \sum\limits^{k}_{j=1}\sum\limits^{r_{j}}_{i=1}\sum\limits^{+\infty}_{s=0}(A t)^{s}\sum\limits^{m_{ji}}_{l=1} v_{l}^{ji}r_{l}
%  \\
 % &=&
 =
  \sum\limits^{k}_{j=1}\sum\limits^{r_{j}}_{i=1}\sum\limits^{m_{ji}}_{l=1}\sum\limits^{+\infty}_{s=0}(A t)^{s}v_{l}^{ji}r_{l},
  \end{array}
\label{solve1}
\end{equation}
where $\eta\in   \mathbb{H}^{n} $. According to Eq.(\ref{vl}), then
\begin{equation}
\begin{array}{lll}
(A t)^{0}v_{l}^{ji}r_{l}&=&v_{l}^{ji}r_{l},
\\
(A t)^{1}v_{l}^{ji}r_{l}&=&t(v_{l-1}^{ji}+v_{l}^{ji}\lambda_{j})r_{l},
\\
(A t)^{2}v_{l}^{ji}r_{l}&=&\frac{t^{2}}{2!}(v_{l-2}^{ji}+2v_{l-1}^{ji}\lambda_{j}+v_{l}^{ji}\lambda_{j}^{2})r_{l},
\\
(At)^{3}v_{l}^{ji}r_{l}&=&\frac{t^{3}}{3!}(C_{3}^{3}v_{l-3}^{ji}+C_{3}^{2}v_{l-2}^{ji}\lambda_{j}+
C_{3}^{1}v_{l-1}^{ji}\lambda_{j}^{2}+v_{l}^{ji}\lambda_{j}^{3})r_{l},
\\
&\,\,\vdots &
\\
(At)^{l-1}v_{l}^{ji}r_{l}&=&\frac{t^{l-1}}{(l-1)!}(C_{l-1}^{l-1}v_{1}^{ji}+C_{l-1}^{l-2}v_{2}^{ji}\lambda_{j}+\cdots+C_{l-1}^{1}v_{l-1}^{ji}\lambda_{j}^{l-2}+v_{p}^{ji}\lambda_{j}^{l-1})r_{l},
\\
(At)^{l}v_{l}^{ji}r_{l}&=&\frac{t^{l}}{(l)!}(C_{l}^{l-1}v_{1}^{ji}\lambda_{j}+C_{l}^{l-2}v_{2}^{ji}\lambda_{j}^{2}+\cdots+C_{l}^{1}v_{l-1}^{ji}\lambda_{j}^{l-1}+v_{p}^{ji}\lambda_{j}^{l})r_{l},
\\
(At)^{l+1}v_{l}^{ji}r_{l}&=&\frac{t^{l+1}}{(l+1)!}(C_{l+1}^{l-1}v_{1}^{ji}\lambda_{j}^{2}+C_{l+1}^{l-2}v_{2}^{ji}\lambda_{j}^{3}+\cdots+C_{l+1}^{1}v_{l-1}^{ji}\lambda_{j}^{l}+v_{p}^{ji}\lambda_{j}^{l+})r_{l},
\\
&\,\,\vdots &
\end{array}
\label{solve2}
\end{equation}
Substituting Eq.(\ref{solve2}) into $\sum\limits^{+\infty}_{s=0}(A t)^{s}v_{l}^{ji}r_{l} $, we can get
\begin{equation}
\begin{array}{lll}
\sum\limits^{+\infty}_{s=0}(A t)^{s}v_{l}^{ji}r_{l}=(v_{l}^{ji}+tv_{l-1}^{ji}+\frac{t^{2}}{2!}v_{l-2}^{ji}+\cdots+\frac{t^{l-2}}{(l-2)!}v_{2}^{ji}+
\frac{t^{l-1}}{(l-1)!}v_{1}^{ji})(\exp\,\lambda_{j} t)r_{l}.
\end{array}
\label{solve3}
\end{equation}
Consequently, we have
\begin{equation}
\begin{array}{lll}
x(t)=(\exp\,At) \eta=\sum\limits^{k}_{j=1}\sum\limits^{r_{j}}_{i=1}\sum\limits^{m_{ji}}_{l=1}(v_{l}^{ji}+tv_{l-1}^{ji}+\frac{t^{2}}{2!}v_{l-2}^{ji}+
\cdots+\frac{t^{l-1}}{(l-1)!}v_{1}^{ji})(\exp\,\lambda_{j} t)r_{l}.
\end{array}
\label{solve4}
\end{equation}

Therefore, if the $A$ is real or complex matrix, the form of solution $x(t)=\exp\{A t\} \eta$ of Eq.(\ref{Constant}) is the same as that ordinary form.

Secondly, how to get the solution $x(t)=\exp\{A t\} \eta$ of Eq.(\ref{solve4})? If we get the eigenvalue $\lambda_{j}$, By $Av_{1}^{ji}=v_{1}^{ji}\lambda_{j}$ and Eq.(\ref{vl}), we can get the set $\{v_{1_{i}}^{j}, v_{2_{i}}^{j},\cdots,v_{m_{ji}-1}^{j},v_{m_{ji}}^{j}\}$. For computational convenience, we introduce a
method to compute the eigenvalue $\lambda_{j}$ and the set $\{v_{1_{i}}^{j}, v_{2_{i}}^{j},\cdots,v_{m_{ji}-1}^{j},v_{m_{ji}}^{j}\}$. First, we will introduce the following results.

For $A\in\mathbb{H}^{n\times n}$, $A=A_{1}+A_{2}\boldsymbol{j}$, where $A_{1}$ and $A_{2}$ are $n\times n$ complex matrices. We associate with $A$ the $2n\times 2n$ complex matrix
\begin{equation}
\phi(A)=\left ( \begin{array}{ll}
\,\,\,\,A_{1}  & A_{2}
 \\
-\overline{A_{2}} &  \overline{A_{1}}
\end{array}
\right ),
\label{AA}
\end{equation}
and call $\phi(A)$ the complex adjoint matrix of the quaternion matrix $A$. Let $\Sigma$ be the collection of all $2n\times 2n$ partitioned complex matrices in the form (\ref{AA}).

For $v\in\mathbb{H}^{n}$, write $v=v_{1}+v_{2}\boldsymbol{j}$, where $v_{1}$ and $v_{2}$ are complex n-tuples . We associate with $v$ the complex $2n$-tuples
\begin{equation}
\varphi(v)=\left ( \begin{array}{ll}
\,\,\,\,v_{1}
 \\
-\overline{v_{2}}
\end{array}
\right ).
\label{vv}
\end{equation}

The mapping $v\rightarrow\varphi(v)$ is an isomorphism between $\mathbb{H}^{n}$ and $\mathbb{C}^{2n}$ obviously. And $v\neq 0$ if and only if $\varphi(v)\neq 0$ can be easily obtained. And $\varphi(v)^{\ast}$ is called the adjoint vector of the $\varphi(v)$
\[
\varphi(v)^{\ast}=
\left ( \begin{array}{ll}
\,\,\,\,v_{1}
 \\
-\overline{v_{2}}
\end{array}
\right )^{\ast}
=\left ( \begin{array}{ll}
v_{2}
 \\
\overline{v_{1}}
\end{array}
\right ).
\label{vv2}
\]

\begin{lemma}\label{Le6.7}
For $A\in\mathbb{H}^{n\times n}$, $v,u\in\mathbb{H}^{n}$ and $\lambda\in\mathbb{C}$, if $\phi(A)\varphi(v)=\varphi(u)+\varphi(v)\lambda$ holds, then $\phi(A)\varphi(v)^{\ast}=\varphi(u)^{\ast}+\varphi(v)^{\ast}\overline{\lambda}$ holds.
\end{lemma}

It can be easily proved by Lemma 3 \cite{Zhangf}.

\begin{lemma}\label{Le6.8}
For $A\in\mathbb{H}^{n\times n}$, $v,u\in\mathbb{H}^{n}$ and $\lambda\in\mathbb{C}$, if $Av=u+v\lambda$ holds, if and only if $\phi(A)\varphi(v)=\varphi(u)+\varphi(v)\lambda$ holds.
\end{lemma}

{\bf Proof.} Note that
\[
Av=(A_{1}+A_{2}\boldsymbol{j})(v_{1}+v_{2}\boldsymbol{j})
=A_{1}v_{1}+A_{2}\boldsymbol{j}v_{2}\boldsymbol{j}+A_{2}\boldsymbol{j}v_{1}+A_{1}v_{2}\boldsymbol{j},
\]
\[
u+v\lambda=(u_{1}+u_{2}\boldsymbol{j})+(v_{1}+v_{2}\boldsymbol{j})\lambda=
u_{1}+v_{1}\lambda+(u_{2}\boldsymbol{j}+v_{2}\boldsymbol{j}\lambda).
\]
If $Av=u+v\lambda$ holds, then
\begin{equation}
\begin{array}{ccc}
A_{1}v_{1}+A_{2}\boldsymbol{j}v_{2}\boldsymbol{j}&=&u_{1}+v_{1}\lambda,
\\
A_{2}\boldsymbol{j}v_{1}+A_{1}v_{2}\boldsymbol{j}&=&u_{2}\boldsymbol{j}+v_{2}\boldsymbol{j}\lambda,
\label{Avu1}
\end{array}
\end{equation}
which implies
\begin{equation}
\begin{array}{ccc}
A_{1}v_{1}-A_{2}\overline{v_{2}}&=&u_{1}+v_{1}\lambda,
\\
\overline{A_{2}}v_{1}+\overline{A_{1}}\overline{v_{2}}&=&\overline{u_{2}}+\overline{v_{2}}\lambda.
\label{Avu2}
\end{array}
\end{equation}
It follows that
\begin{equation}
\left ( \begin{array}{ll}
\,\,\,\,A_{1}  & A_{2}
 \\
-\overline{A_{2}} &  \overline{A_{1}}
\end{array}
\right )
\left ( \begin{array}{ll}
\,\,\,\,v_{1}
 \\
-\overline{v_{2}}
\end{array}
\right )=
\left ( \begin{array}{ll}
\,\,\,\,u_{1}
 \\
-\overline{u_{2}}
\end{array}
\right )+
\left ( \begin{array}{ll}
\,\,\,\,v_{1}
 \\
-\overline{v_{2}}
\end{array}
\right )\lambda.
\label{Avu3}
\end{equation}
Therefore,
\[
\phi(A)\varphi(v)=\varphi(u)+\varphi(v)\lambda.
\]
Conversely, it can be easily proved.
\begin{corollary}\label{Co6.3}
For $A\in\mathbb{H}^{n\times n}$, $v\in\mathbb{H}^{n}$ and $\lambda\in\mathbb{C}$, $Av=v\lambda$ holds, if and only if $\phi(A)\varphi(v)=\varphi(v)\lambda$ holds.
\end{corollary}

To obtain $\lambda_{j}$ and the set $\{v_{1_{i}}^{j}, v_{2_{i}}^{j},\cdots,v_{m_{ji}-1}^{j},v_{m_{ji}}^{j}\}$, we introduce the computational process of this method according to the proof of Theorem 1 \cite{Zhangf}.
Let $\lambda_{j}=a+b\boldsymbol{i}$ (k-fold) is a eigenvalue of $\phi(A)$, $A\in\mathbb{H}^{n\times n}$. By Corollary \ref{Co6.3} and Theorem 1 \cite{Zhangf}. Then
\begin{description}
\item[(i)] If $b>0$, $\lambda_{j}=a+b\boldsymbol{i}$ (k-fold) is a eigenvalue of $A$ and set $\{\varphi(v_{1_{i}}^{j}),\varphi(v_{2_{i}}^{j}),\cdots,\varphi(v_{m_{ji}-1}^{j}),
    \varphi(v_{m_{ji}}^{j})\}$ is easily calculated. (the set $\{v_{1_{i}}^{j}, v_{2_{i}}^{j},\cdots,v_{m_{ji}-1}^{j},v_{m_{ji}}^{j}\}$ is undetermined).  By Eq.(\ref{vv}), we can obtain the set $\{v_{1_{i}}^{j}, v_{2_{i}}^{j},\cdots,v_{m_{ji}-1}^{j},v_{m_{ji}}^{j}\}$.
\item[(ii)] If $b=0$, $\lambda_{j}=a+b\boldsymbol{i}$ ($\frac{k}{2}$-fold) is a eigenvalue of $A$, and by Lemma \ref{Le6.6}, there two sets $\{\varphi(v_{1_{i}}^{j}), \varphi(v_{2_{i}}^{j}),\cdots,\varphi(v_{m_{ji}-1}^{j}),\varphi(v_{m_{ji}}^{j})\}$ and $\{\varphi(v_{1_{i}}^{j})^{\ast},\varphi(v_{2_{i}}^{j})^{\ast},\cdots,\varphi(v_{m_{ji}-1}^{j})^{\ast},
    \varphi(v_{m_{ji}}^{j})^{\ast}\}$ are calculated. Taking the set   $\{\varphi(v_{1_{i}}^{j}), \varphi(v_{2_{i}}^{j})  \cdots,\varphi(v_{m_{ji}-1}^{j}),\varphi(v_{m_{ji}}^{j})\}$, By Eq.(\ref{vv}), we can obtain the set $\{v_{1_{i}}^{j}, v_{2_{i}}^{j}, \cdots,v_{m_{ji}-1}^{j},v_{m_{ji}}^{j}\}$.
\end{description}

Now we are in a position to obtain $exp\{At\}$ from Eq.(\ref{solve4}), we can firstly choose $n$ independent initial value vector, then the corresponding $n$ solutions to the IVP are independent. For convenience, we usually choose the natural basis. Let $\eta=e_{1},\eta=e_{2},\cdots,\eta=e_{n}$, correspondingly, we can get $n$ independent solutions. These $n$ independent solutions compose the column of $exp\{At\}$.
Noticing that $exp\{At\}=exp\{At\}E=[(exp\{At\})e_{1},(exp\{At\})e_{2},\cdots,(exp\{At\})e_{n}]$,
where
\[
e_{1}=\left [\begin{array}{l} 1
\\
0
\\
\vdots
\\
0
\\
0
\end{array}
\right ],
e_{2}=\left [\begin{array}{l} 0
\\
1
\\
\vdots
\\
0
\\
0
\end{array}
\right ],
\cdots,
e_{n}=\left [\begin{array}{l} 0
\\
0
\\
\vdots
\\
0
\\
1
\end{array}
\right ],
\]
are the unit vector.

Some examples are presented to show the validity of this method.

\noindent {\bf Example 6.2} Find a fundamental matrix of the following QDEs
\begin{equation}
\dot {x}= \left ( \begin{array}{ll}
\mathbf{ i}  &  1
 \\
 0  &  \mathbf{j}
\end{array}
\right )x,\,\,\,\, x=(x_1,x_2)^T.
\label{Ex6.3.3}
\end{equation}

{\bf Answer:} From Lemma 6.8 and Lemma 6.9 in \cite{QDE1}, we see that $\lambda_1=\mathbf{i}$ and $\lambda_2=\mathbf{j}$. It should be noted that $\mathbf{j}$ is similar to $\mathbf{i}$. In fact, taking $\alpha=1+\mathbf{i}+\mathbf{j}+\mathbf{k}$, then $\alpha^{-1}=\frac{1}{4}(1-\mathbf{i}-\mathbf{j}-\mathbf{k})$. Consequently, $\mathbf{j}=\alpha^{-1}\mathbf{i}\alpha$, that is, $\mathbf{j}$ is similar to $\mathbf{i}$.
To find the eigenvector of $\lambda_1=i$, we consider the following equation
\[
Aq=q\lambda_1,
\]
that is
\begin{equation}
 \left\{ \begin{array}{ccc}
 \mathbf{i}  q_1 + q_2 & =& q_1 \mathbf{i}
 \\
\mathbf{j} q_2 &=& q_2 \mathbf{i}.
\end{array}
\right.
\label{EV1}
\end{equation}

Let $q_{1}=a_{0}+a_{1}\boldsymbol{i}+a_{2}\boldsymbol{j}+a_{3}\boldsymbol{k}$ and $q_{2}=b_{0}+b_{1}\boldsymbol{i}+b_{2}\boldsymbol{j}+b_{3}\boldsymbol{k}$. Substituting it into the Eq.(\ref{EV1}), we obtain $q_1=a_{0}+a_{1}\boldsymbol{i},\,q_{2}=0\,\,a_{0},a_{1}\in \mathbb{R}$.

According to Remark \ref{Re6.1.1}, it is impossible to find two independent eigenvector of (\ref{Ex6.3.3}).

We consider matrix $\phi(A)$, the eigenvalue of $\phi(A)$ are $\lambda_{1}=\boldsymbol{i}$ (2-fold), $\lambda_{2}=-\boldsymbol{i}$ (2-fold). The eigenvector of $\lambda_{1}=\boldsymbol{i}$ is $\varphi(v)=(1,0,0,0)^{T}$, and
\[
v=\left ( \begin{array}{c}
1
 \\
0
\end{array}
\right ).
\]
By Lemma \ref{Le6.8} and $\phi(A)\varphi(u)=\varphi(v)+\varphi(u)\lambda_{1}$, we can get $\varphi(u)=(0,1,\frac{1}{2},\boldsymbol{i})^{T}$ and
\[
u=\left ( \begin{array}{c}
-\frac{1}{2}\boldsymbol{j}
 \\
1+\boldsymbol{k}
\end{array}
\right ).
\]
Substituting $v,u$ into Eq.(\ref{solve4}), For any solutions $exp\{At\}\eta$, let $\eta=v r_1+u r_2$ we can get
\[
exp\{At\}\eta=v e^{\boldsymbol{i} t} r_1+(u+v t)e^{\boldsymbol{i} t} r_2.
\]
Namely
\[
exp\{At\}\eta=
\left ( \begin{array}{ccc}
1 & -\frac{1}{2}\boldsymbol{j}+t
 \\
0 &1+\boldsymbol{k}
\end{array}
\right )
\left ( \begin{array}{c}
 e^{\boldsymbol{i} t} r_1
 \\
 e^{\boldsymbol{i} t} r_2
\end{array}
\right ).
\]

Let $\eta=(1, 0)^T,(0, 1)^T$ in turn, we can obtain two linear
independent solutions, which compose the fundamental matrix $exp\{At\}$, namely
\[
exp\{At\}=
\left ( \begin{array}{ccc}
e^{\boldsymbol{i}t} & e^{\boldsymbol{i}t}\frac{1-\boldsymbol{k}}{2}-(\frac{1}{2}\boldsymbol{j}-t)e^{\boldsymbol{i}t}\frac{-\boldsymbol{i}+\boldsymbol{j}}{4}
 \\
0 &(1+\boldsymbol{k})e^{\boldsymbol{i}t}\frac{-\boldsymbol{i}+\boldsymbol{j}}{4}
\end{array}
\right ).
\]

\noindent {\bf Example 6.3} Find fundamental matrix $exp\{At\}$ of the following QDEs
\begin{equation}
\dot {x}=Ax= \left ( \begin{array}{ccc}
 \boldsymbol{i} & \boldsymbol{j} & \boldsymbol{j}
 \\
 \boldsymbol{k} & 1 & \boldsymbol{k}
 \\
 0 & 0 & 1
\end{array}
\right )x,\,\,\,\, x=(x_1,x_2,x_3)^T.
\label{Ex6.4}
\end{equation}

{\bf Answer:} we can easily get the eigenvalue of $\phi(A)$ are $\lambda_{1}=0$ (2-fold), $\lambda_{2}=1+\boldsymbol{i}$, $\lambda_{3}=1-\boldsymbol{i}$ and $\lambda_{4}=1$ (2-fold). The eigenvector of $\lambda_{1}=0$ is $\varphi(v_{1})=(0,-\boldsymbol{i},0,1,0,0)^{T}$. By Eq.(\ref{vv}), we get
\[
v_{1}=\left ( \begin{array}{c}
 -\boldsymbol{j}
 \\
-\boldsymbol{i}
 \\
0
\end{array}
\right ).
\]
The eigenvectors of $\lambda_{2}=1+\boldsymbol{i}$ are $\varphi(v_{2})=(1,0,0,0,1,0)^{T}$ and
\[
v_{2}=\left ( \begin{array}{c}
 1
 \\
-\boldsymbol{j}
 \\
0
\end{array}
\right ).
\]
The eigenvectors of $\lambda_{4}=1$ are $\varphi(v_{3})=(\frac{1+\boldsymbol{i}}{2},\frac{1+\boldsymbol{i}}{2},-\frac{1+\boldsymbol{i}}{2},0,1,0)^{T}$ and
\[
v_{3}=\left ( \begin{array}{c}
 \frac{1+\boldsymbol{i}}{2}
 \\
\frac{1+\boldsymbol{i}}{2}-\boldsymbol{j}
 \\
-\frac{1+\boldsymbol{i}}{2}
\end{array}
\right ).
\]
%Substituting $v_{1},v_{2},v_{2}$ into Eq.(\ref{solve4}), For any solutions $exp\{At\}\eta$, let $\eta=v_{1} r_1+v_{2} r_2+u r_3$ we can get
%\[
%exp\{At\}\eta=v_{1} e^{\boldsymbol{i} t} r_1+v_{2} e^{t} r_2+v_{3} e^{t} r_3.
%\]
%Namely
%\[
%exp\{At\}\eta=
%\left ( \begin{array}{ccc}
%\frac{1+\boldsymbol{i}}{2} & -\boldsymbol{j} & 1
 %\\
% \frac{1+\boldsymbol{i}}{2}-\boldsymbol{j} &-\boldsymbol{i} & -\boldsymbol{j}
% \\
%- \frac{1+\boldsymbol{i}}{2} & 0 & 0
%\end{array}
%\right )
%\left ( \begin{array}{c}
% e^{t} r_1
 %\\
 %r_2
 %\\
 %e^{(1+\boldsymbol{i})t} r_3
%\end{array}
%\right ).
%\]
%Set $\eta=v_{1},\eta=v_{2},\eta=v_{3}$. Since $v_{1},v_{2},v_{3}$ are right independent, we obtain three independent solutions, which compose the fundamental matrix.
%\[
%\left ( \begin{array}{ccc}
%\frac{1+\boldsymbol{i}}{2}e^{t} & -\boldsymbol{j} & e^{(1+\boldsymbol{i})t}
% \\
% (\frac{1+\boldsymbol{i}}{2}-\boldsymbol{j})e^{t} &-\boldsymbol{i} & -\boldsymbol{j}e^{(1+\boldsymbol{i})t}
% \\
%- \frac{1+\boldsymbol{i}}{2}e^{t} & 0 & 0
%\end{array}
%\right ).
%\]
From corollary \ref{Co6.1}, then
\[
M(t)=\left ( \begin{array}{ccc}
-\boldsymbol{j} & e^{(1+\boldsymbol{i})t} & \frac{1+\boldsymbol{i}}{2}e^{t}
 \\
-\boldsymbol{i} & -\boldsymbol{j}e^{(1+\boldsymbol{i})t} &  (\frac{1+\boldsymbol{i}}{2}-\boldsymbol{j})e^{t}
 \\
0 & 0 & - \frac{1+\boldsymbol{i}}{2}e^{t}
\end{array}
\right ).
\]
is a fundamental matrix of Eq.(\ref{Ex6.4}).

For any solutions $exp\{At\}\eta$ of Eq.(\ref{Ex6.4}). Substituting $v_{1},v_{2},v_{3}$ into Eq.(\ref{solve4}), let $\eta=v_{1} r_1+v_{2} r_2+v_{3} r_3$ we can get
\[
exp\{At\}\eta=v_{1} r_1+v_{2} e^{(1+\boldsymbol{i})t} r_2+v_{3} e^{t} r_3,
\]
namely
\[
exp\{At\}\eta=
\left ( \begin{array}{ccc}
-\boldsymbol{j} & 1 & \frac{1+\boldsymbol{i}}{2}
 \\
-\boldsymbol{i} & -\boldsymbol{j} &  \frac{1+\boldsymbol{i}}{2}-\boldsymbol{j}
 \\
0 & 0 & - \frac{1+\boldsymbol{i}}{2}
\end{array}
\right )
\left ( \begin{array}{c}
r_1
 \\
 e^{(1+\boldsymbol{i})t} r_2
 \\
 e^{t} r_3
\end{array}
\right ).
\]

Let $\eta=(1, 0, 0)^T,(0, 1, 0)^T,(0, 1, 0)^T$ in turn, we can obtain three linear
independent solutions, which compose the fundamental matrix $exp\{At\}$, namely
\[
exp\{At\}=
\left ( \begin{array}{ccc}
  \frac{1-\boldsymbol{i}}{2}+\frac{1+\boldsymbol{i}}{2}\alpha& \frac{-\boldsymbol{j}+\boldsymbol{k}}{2}+\beta&
  \boldsymbol{j}\gamma+\delta-e^{t}
 \\
\frac{\boldsymbol{j}-\boldsymbol{k}}{2}-\frac{\boldsymbol{j}-\boldsymbol{k}}{2}\alpha &\frac{1-\boldsymbol{i}}{2}-\boldsymbol{j}\beta & \boldsymbol{i}\gamma-\boldsymbol{j}\delta-(1-\boldsymbol{j}-\boldsymbol{k})e^{t}
 \\
0 & 0 & e^{t}
\end{array}
\right ).
\]
where $\alpha= e^{(1+\boldsymbol{i})t}$, $\beta= e^{(1+\boldsymbol{i})t}\frac{\boldsymbol{j}-\boldsymbol{k}}{2}$, $\gamma=\frac{1+\boldsymbol{i}+
\boldsymbol{j}-\boldsymbol{k}}{2}$, $\delta= e^{(1+\boldsymbol{i})t}\frac{1-\boldsymbol{i}+\boldsymbol{j}-\boldsymbol{k}}{2}$.

\noindent {\bf Example 6.4} Find fundamental matrix $exp\{At\}$ of the following QDEs
\begin{equation}
\dot {x}= Ax=\left ( \begin{array}{ccc}
 \boldsymbol{j} & \boldsymbol{k} & \boldsymbol{i}
 \\
 0 & 1 & \boldsymbol{k}
 \\
 0 & 0 & 1
\end{array}
\right )x,\,\,\,\, x=(x_1,x_2,x_3)^T.
\label{Ex6.5}
\end{equation}

{\bf Answer:} we easily get the eigenvalue of $\phi(A)$ are $\lambda_{1}=\boldsymbol{i}$, $\lambda_{2}=-\boldsymbol{i}$ and $\lambda_{3}=1$ (4-fold). The eigenvector of $\lambda_{1}=\boldsymbol{i}$ is $\varphi(v_{1})=(-\boldsymbol{i},0,0,1,0,0)^{T}$, and
\[
v_{1}=\left ( \begin{array}{c}
 -\boldsymbol{i}-\boldsymbol{j}
 \\
0
 \\
0
\end{array}
\right ).
\]

The eigenvectors of $\lambda_{3}=1$ are $\varphi(v_{2})=(-\boldsymbol{i},1,0,0,1,0)^{T}$ and $\varphi(v_{3})=\varphi(v_{2})^{\ast}$, and

\[
v_{2}=\left ( \begin{array}{c}
 \boldsymbol{i}
 \\
1-\boldsymbol{j}
 \\
0
\end{array}
\right ).
\]
By Lemma \ref{Le6.8} and $\phi(A)\varphi(u)=\varphi(v_{2})+\varphi(u)\lambda$, we can get $\varphi(u)=(1-\boldsymbol{i},-1-2\boldsymbol{i},-\boldsymbol{i},0,0,-\boldsymbol{i})^{T}$ and
\[
u=\left ( \begin{array}{c}
1-\boldsymbol{i}
 \\
-1-2\boldsymbol{i}
 \\
-\boldsymbol{i}-\boldsymbol{k}
\end{array}
\right ).
\]
Let $\eta=v_{1} r_1+v_{2} r_2+u r_3$. For any solutions $exp\{At\}\eta$ of Eq.(\ref{Ex6.5}), substituting $v_{1},v_{2},u$ into Eq.(\ref{solve4}), we can get
\[
exp\{At\}\eta=v_{1} e^{\boldsymbol{i} t} r_1+v_{2} e^{t} r_2+(u+v_{2} t) e^{t} r_3.
\]
That is,
\[
exp\{At\}\eta=
\left ( \begin{array}{ccc}
 -\boldsymbol{i}-\boldsymbol{j} & \boldsymbol{i} & 1-\boldsymbol{i}+\boldsymbol{i} t
 \\
0 &1-\boldsymbol{j} & 1-2\boldsymbol{i}-(\boldsymbol{i}+\boldsymbol{k}) t
 \\
0 & 0 & -\boldsymbol{i}-\boldsymbol{k}
\end{array}
\right )
\left ( \begin{array}{c}
 e^{\boldsymbol{i} t} r_1
 \\
 e^{t} r_2
 \\
 e^{t} r_3
\end{array}
\right ).
\]
Taking $\eta=(1, 0, 0)^T,(0, 1, 0)^T,(0, 1, 0)^T$ in turn, we can obtain three linear
independent solutions, which compose the fundamental matrix $exp\{At\}$, namely
\[
exp\{At\}=
\left ( \begin{array}{ccc}
 \alpha\frac{\boldsymbol{i}+\boldsymbol{j}}{2} & \alpha\frac{1-\boldsymbol{i}+\boldsymbol{j}+\boldsymbol{k}}{4} +\frac{\boldsymbol{i}+\boldsymbol{k}}{2}e^{t} & \alpha\frac{2+\boldsymbol{i}-\boldsymbol{j}}{2}+ (\frac{\boldsymbol{i}+\boldsymbol{j}-\boldsymbol{k}-t-\boldsymbol{j}t}{2}) e^{t}
 \\
0 &e^{t} & \boldsymbol{k}t e^{t}
 \\
0 & 0 & e^{t}
\end{array}
\right ).
\]
where $\alpha= -(\boldsymbol{i}+\boldsymbol{j})e^{\boldsymbol{i} t}$.

%To obtain the fundamental matrix, set $\eta=v_{1},\eta=v_{2},\eta=u$. Since $v_{1},v_{2},u$ are right independent, we can get three independent solutions, which compose the fundamental matrix.
%\[
%\left ( \begin{array}{ccc}
% -(\boldsymbol{i}+\boldsymbol{j})e^{\boldsymbol{i} t} & \boldsymbol{i} e^{t} & (1-\boldsymbol{i}+\boldsymbol{i} t) e^{t}
% \\
%0 &1-\boldsymbol{j} e^{t} & (1-2\boldsymbol{i}-(\boldsymbol{i}+\boldsymbol{k}) t) e^{t}
% \\
%0 & 0 & -(\boldsymbol{i}+\boldsymbol{k}) e^{t}
%\end{array}
%\right ).
%\]

\section{Conclusion  }

In this paper, we presented an algorithm to evaluate the fundamental matrix by employing the eigenvalue and eigenvectors. We gave a method to construct the fundamental matrix when the linear system has multiple eigenvalues. In particular, if the number of independent eigenvectors might be smaller than the dimensionality of the system. That is, the numbers of the eigenvectors is not enough to construct a fundamental matrix. We therefore have to discover how to find the ``missing solutions". The main purpose is to answer this question.

\smallskip

\section{Conflict of Interests}

The authors declare that there is no conflict of interests
regarding the publication of this article.

\end{document}